\documentclass[final]{siamltex}

\usepackage{epsfig,amssymb,latexsym}
\usepackage{amsfonts,psfrag,amsmath,bbm,color}
\usepackage{cancel}

\usepackage{mathrsfs}
\usepackage{graphicx}
\usepackage{textcomp}
\usepackage{multirow}
\usepackage{enumerate}
\usepackage{cancel}
\usepackage{caption}  
\usepackage{url}

\usepackage{algorithm}
\usepackage{algorithmic}

\usepackage{mathtools}

\DeclarePairedDelimiter\floor{\lfloor}{\rfloor}

\usepackage{amscd}

\graphicspath{{./figs/}}

\newtheorem{remark}{Remark}[section]

\def\PP{{{\rm l}\kern - .15em {\rm P} }}
\def\PN2{{\PP_{N}-\PP_{N-2}}}

\newcommand{\R}{\mathbbm{R}}

\newcommand{\bphi}{\boldsymbol{\varphi}}

\newcommand{\btau}{\boldsymbol{\tau}}

\newcommand{\ba}{\boldsymbol{a}}

\newcommand{\bu}{\boldsymbol{u}}

\newcommand{\bur}{{\boldsymbol{u}}_r}

\newcommand{\bx}{\boldsymbol{x}}
\newcommand{\bX}{\boldsymbol{X}}

\newcommand{\bXr}{{\bf X}^r}

\newcommand{\tA}{\tilde{A}}
\newcommand{\tB}{\widetilde{B}}

\newcommand{\obu}{\overline{\boldsymbol u}}

\newcommand{\blue}[1]{{\color{blue}#1}}

\newcommand{\deleted}[1]{{}}

\begin{document}

\title{Physically-Constrained Data-Driven, Filtered Reduced Order Modeling of Fluid Flows}

\author{
M. Mohebujjaman
\thanks{Department of Mathematics, Virginia Tech, Blacksburg, VA, 24061; email: jaman@vt.edu}
\and L. G. Rebholz
\thanks{Department of Mathematical Sciences, Clemson University, Clemson, SC 29634; 
Partially supported by NSF DMS1522191 and Army Research Office 65294-MA, email: rebholz@clemson.edu}
\and T. Iliescu
\thanks{Department of Mathematics, Virginia Tech, Blacksburg, VA 24061
Partially supported by NSF DMS1522656, email: iliescu@vt.edu}.
}

\date{\today}

\maketitle

\begin{abstract}
In~\cite{xie2018data}, we proposed a data-driven filtered reduced order model (DDF-ROM) framework for the numerical simulation of fluid flows, which can be formally written as
\begin{equation*}
	\boxed{
		\text{
			DDF-ROM = Galerkin-ROM + Correction
		}
	}
\end{equation*}
The new DDF-ROM was constructed by using ROM spatial filtering and data-driven ROM closure modeling (for the Correction term) and was successfully tested in the numerical simulation of a 2D channel flow past a circular cylinder at Reynolds numbers $Re=100, Re=500$ and $Re=1000$.

In this paper, we propose a {\it physically-constrained} DDF-ROM (CDDF-ROM), which aims at improving the physical accuracy of the DDF-ROM.
The new physical constraints require that the CDDF-ROM operators satisfy the same type of physical laws (i.e., the nonlinear operator should conserve energy and the ROM closure term should be dissipative) as those satisfied by the fluid flow equations.
To implement these physical constraints, in the data-driven modeling step of the DDF-ROM, we replace the unconstrained least squares problem with a constrained least squares problem. 
We perform a numerical investigation of the new CDDF-ROM and standard DDF-ROM for a 2D channel flow past a circular cylinder at Reynolds numbers $Re=100, Re=500$ and $Re=1000$.
To this end, we consider a reproductive regime as well as a predictive (i.e., cross-validation) regime in which we use as little as $50\%$ of the original training data. 
The numerical investigation clearly shows that the new CDDF-ROM is significantly more accurate than the DDF-ROM in both regimes.
\end{abstract}

\begin{keywords} 
reduced order modeling, 
data-driven modeling, 
spatial filter,
physical constraints
\end{keywords}


\begin{AMS}
65M60, 76F65
\end{AMS}

\pagestyle{myheadings}
\thispagestyle{plain}
\markboth{M. MOHEBUJJAMAN, L. G. REBHOLZ, AND T. ILIESCU}{Physically-Constrained Data-Driven, Filtered Reduced Order Modeling of Fluid Flows}

\section{Introduction}

To present the new reduced order model (ROM), we use the incompressible Navier-Stokes equations (NSE):
\begin{eqnarray}
    && \frac{\partial \bu}{\partial t}
    - Re^{-1} \Delta \bu
    + \bu \cdot \nabla \bu
    + \nabla p
    = {\bf 0} \, ,
    \label{eqn:nse-1}                                                         \\
    && \nabla \cdot \bu
    = 0 \, ,
    \label{eqn:nse-2}
\end{eqnarray}
where $\bu$ is the velocity, $p$ the pressure, and $Re$ the Reynolds number. 
We use the initial condition $\bu(\bx, 0) = \bu_0(\bx)$ and (for simplicity) homogeneous Dirichlet boundary conditions: $\bu(\bx, t) = {\bf 0}$.
ROMs have been used to reduce the computational cost of scientific and engineering applications that are governed by relatively few recurrent dominant spatial structures~\cite{ballarin2016fast,benner2015survey,bistrian2015improved,gunzburger2017ensemble,hesthaven2015certified,HLB96,noack2011reduced,perotto2017higamod,quarteroni2015reduced,stefanescu2015pod}.
In an offline stage, full order model (FOM) data on a given time interval are used to build the ROM.
In an online stage, ROMs are repeatedly used for parameter settings and/or time intervals that are different from those used to build them. 
\\[-0.3cm]

For a given general PDE, {\it Projection ROMs (Proj-ROMs)}~\cite{ballarin2015supremizer,HLB96,noack2011reduced} strategy for approximating its solution, $\bu$, is:
(i) Choose a few dominant modes $\{ \bphi_1, \ldots, \bphi_r \}$ (which represent the recurrent spatial structures) as  basis functions.
(ii) Replace $\bu$ with $\bur = \sum_{j=1}^{r} a_j \, \bphi_j$ in the given PDE.
(iii) Use a Galerkin projection of PDE($\bur$) onto the ROM space $\text{span} \{ \bphi_1, \ldots, \bphi_r \}$ to obtain the Proj-ROM. 
For example, in fluid dynamics, the Proj-ROM often takes the form
\begin{eqnarray}
	\dot{\ba}
	= A \, \ba
	+ \ba^{\top} \, B \, \ba \, ,
	\label{eqn:proj-rom}
\end{eqnarray}
where $\ba$ is the vector of unknown ROM coefficients and $A \in \R^{r \times r}, B \in \R^{r \times r \times r}$ are ROM operators that are assembled in the offline stage. \\[-0.3cm]

\noindent{\it Data-Driven ROMs (DD-ROMs)} (e.g., sparse identification of nonlinear dynamics~\cite{brunton2016discovering} and operator inference method~\cite{peherstorfer2015dynamic,peherstorfer2016data}) use a fundamentally different strategy:
They first postulate a ROM ansatz
\begin{eqnarray}
	\dot{\ba}
	= \tA \, \ba
	+ \ba^{\top} \, \tB \, \ba \, ,
	\label{eqn:dd-rom}
\end{eqnarray}
and then they choose the operators $\tA$ and $\tB$ to minimize the difference between the FOM and ansatz~\eqref{eqn:dd-rom}~\cite{cacuci2013computational,kutz2013data}:
\begin{eqnarray}
		\min_{\widehat{A} , \widehat{B}} \| FOM - ( \dot{\ba} - \widehat{A} \, \ba - \ba^{\top} \, \widehat{B} \, \ba) \|^2 \, .
	\label{eqn:dd-rom-optimization}
\end{eqnarray}
\vspace*{-0.1cm}

\noindent Both Proj-ROMs and DD-ROMs are facing {\it grand challenges}: 
One of the main roadblocks for Proj-ROMs is that they are not accurate models for the dominant modes:
In practice, a corrected Proj-ROM is generally used instead~\cite{amsallem2012stabilization,balajewicz2013low,barone2009stable,benosman2017learning,bergmann2009enablers,carlberg2017galerkin,giere2015supg,osth2014need,protas2015optimal}:
	\begin{eqnarray}
		\dot{\ba}
		= A \, \ba
		+ \ba^{\top} \, B \, \ba
		+ \text{Correction.}
		\label{eqn:proj-rom-modified}
	\end{eqnarray}
Thus, the ROM closure problem (i.e., the modeling of the Correction term in~\eqref{eqn:proj-rom-modified}) needs to be addressed. 
DD-ROMs, on the other hand, can be sensitive to noise in the data, since its operators $\tA$ and $\tB$ are obtained from an inverse problem~\cite{hansen2010discrete,vogel2002computational}.

\vspace*{0.3cm}

In~\cite{xie2018data}, we proposed a {\it data-driven filtered ROM (DDF-ROM)}, which is a {\it hybrid Proj-ROM/DD-ROM framework} in which the classical {\it Proj-ROM} framework is used to model the {\it linear operators} and the {\it DD-ROM} framework is used to model the {\it nonlinear operators}.
Next, we briefly describe the main steps used in the construction of the DDF-ROM.
{\it The philosophy employed in the DDF-ROM construction is to use the classical Galerkin method whenever possible (i.e., for the linear operators) and invoke data-driven modeling only when necessary (i.e., for the nonlinear operators).}
The hybrid character of the new framework is reminiscent of {\it data assimilation}~\cite{kalnay2003atmospheric}: the classical Galerkin method is at the core of the framework, and data-driven modeling is used solely to improve its accuracy.

We build the DDF-ROM framework in two steps.
In the first step, we put forth ROM spatial filtering to discover the {\it exact mathematical formula} for the Correction term in~\eqref{eqn:proj-rom-modified}:
\begin{equation}
	\boxed{
		\text{
			Correction = Exact Mathematical Formula.
		}
	}
	\label{eqn:ddf-rom-step-1}
\end{equation}
In the second step, we utilize data-driven modeling to find a useful approximation for the Correction term in~\eqref{eqn:ddf-rom-step-1}. 
Specifically, we make the ansatz
\begin{eqnarray}
	\text{Correction}
	\approx 
	\tA \, \ba
	+ 
	\ba^{\top} \, \tB \, \ba \, ,
	\label{eqn:ddf-rom-step-2}
\end{eqnarray}
and choose $\tA$ and $\tB$ to minimize the difference between the Exact Mathematical Formula in~\eqref{eqn:ddf-rom-step-1} and our ansatz~\eqref{eqn:ddf-rom-step-2}, both calculated with the available FOM data:
\begin{eqnarray}
	\hspace*{-0.6cm}
	\boxed{
		\min_{\tA , \tB} \| \text{Exact Mathematical Formula}(FOM) - ( \tA \, \ba_{FOM} + \ba_{FOM}^{\top} \, \tB \, \ba_{FOM}) \|^2 \, .
	}
	\label{eqn:f-rom-optimization}
\end{eqnarray}
At the end of the two steps, we obtain the {\it DDF-ROM}
\begin{eqnarray}
	\boxed{
		\dot{\ba}
		= (A + \tA) \, \ba
		+ \ba^{\top} \, (B + \tB) \, \ba \, .
	}
	\label{eqn:ddf-rom-intro}
\end{eqnarray}

The DDF-ROM framework solves the Proj-ROM's closure problem, since the Correction term is modeled by data-driven modeling in~\eqref{eqn:ddf-rom-step-2}--\eqref{eqn:f-rom-optimization}. 
Furthermore, DDF-ROM is more robust to noise than standard DD-ROMs, since DDF-ROM employs data-driven modeling (which is an inverse problem sensitive to noise) only to  model the Correction term in~\eqref{eqn:proj-rom-modified}, whereas the DD-ROMs use it to model {\it all} the ROM operators (compare~\eqref{eqn:f-rom-optimization} with~\eqref{eqn:dd-rom-optimization}).

We note that data-driven closure modeling for non-ROM settings is an extremely active research area, see, e.g., \cite{duraisamy2018turbulence, ling2016machine}.
We also note that there are other data-driven ROM closure models, 
see, e.g., \cite{chekroun2017markovian, kondrashov2015data, lu2017data, gouasmi2017priori, san2018neural, baiges2015reduced, oberai2016approximate}. 
We emphasize, however, that these data-driven ROM closure models are different from our DDF-ROM: 

(i) They do not use spatial filtering (as in LES) to isolate the ROM subfilter-scale stress tensor $\btau$ (i.e., the Correction term);

(ii) They generally do closure modeling for both linear and nonlinear terms \cite{san2018neural, baiges2015reduced, oberai2016approximate}; and 

(iii) They only use the ROM projection to define the Correction term. 
In contrast, our DDF-ROM framework is general and can accommodate {\it any} type of spatial filter.
For example, we could use the ROM differential filter, which was shown to outperform the ROM projection in the numerical investigation of ROMs for 3D flow past a cylinder at $Re=1000$~\cite{wells2017evolve}. 

In~\cite{xie2018data}, we investigated the DDF-ROM in the numerical simulation of a 2D channel flow past a circular cylinder at Reynolds numbers $Re=100, Re=500$, and $Re=1000$.
The DDF-ROM was significantly more accurate than the standard projection ROM.
Furthermore, the computational costs of the DDF-ROM and the standard projection ROM were similar, both costs being orders of magnitude lower than the computational cost of the full order model. 
\\[-0.3cm]

Although the DDF-ROM yielded good results~\cite{xie2018data}, these results got worse when the ROM dimension was  decreased below a certain threshold.
To address this issue, we propose the {\it physically constrained DDF-ROM (CDDF-ROM)}, which aims at improving the physical accuracy of the DDF-ROM.  
These physical constraints require that the operator $\tA$ in the DDF-ROM~\eqref{eqn:ddf-rom-intro} be {\it negative semidefinite} and the operator $\tB$ be {\it energy conserving} (which resembles the constraints satisfied by the operators $A$ and $B$).
To implement these physical constraints, we replace the unconstrained least squares problem~\eqref{eqn:f-rom-optimization} solved in the data-driven modeling step of the DDF-ROM with a constrained least squares problem.

We note that it has long been known that improving the physical accuracy of a discretization leads to more accurate solutions in all measures, especially over long time intervals.  
Thus, physical constraints have been used for decades in the CFD community (see, e.g., Arakawa's pioneering work~\cite{A66}, as well as more recent developments~\cite{AM03,CLRW10,F75,LW04,R07, SCN15,ST89}).
More recently, physical constraints have started to be used in standard (i.e., without ROM) LES closure modeling, see, e.g.,~\cite{duraisamy2018turbulence}.
Finally, physical constraints have also been used in standard ROM (i.e., without closure modeling), see, e.g.,~\cite{loiseau2018constrained, kalashnikova2010stability, rowley2004model}.
To our knowledge, the CDDF-ROM proposed in this paper is the {\it first physically constrained ROM closure model}.
\\[-0.3cm]

The rest of the paper is organized as follows:
In Section~\ref{sec:ddf-rom}, we review the DDF-ROM proposed in~\cite{xie2018data}.
In Section~\ref{sec:cddf-rom}, we propose the new CDDF-ROM.
In Section~\ref{sec:numerical-results}, we perform a numerical investigation of the new CDDF-ROM in the numerical simulation of a 2D flow past a circular cylinder at Reynolds numbers $Re=100, Re=500$, and $Re=1000$.
Finally, in Section~\ref{sec:conclusions}, we summarize our findings and outline future research directions.


\section{Data-Driven Filtered ROM (DDF-ROM)}
	\label{sec:ddf-rom}

In this section, we briefly review the DDF-ROM proposed in~\cite{xie2018data}.
For more details, the reader is referred to~\cite{xie2018data}.

For the ROM basis, we use the proper orthogonal decomposition (POD)~\cite{HLB96,noack2011reduced,Sir87abc,volkwein2013proper}.
We emphasize, however, that other ROM bases (e.g., the dynamic mode decomposition (DMD)~\cite{kutz2016dynamic,mezic2005spectral,rowley2009spectral,schmid2010dynamic}) could be used in the DDF-ROM construction.
Given snapshots (e.g., the finite element solutions) of the NSE~\eqref{eqn:nse-1}--\eqref{eqn:nse-2}, the POD space $\bXr = \text{span} \{ \bphi_{1}, \ldots, \bphi_{r} \}$ approximates the snapshots optimally with respect to the  $L^2$-norm.
The POD approximation of the velocity is defined as 
\begin{equation}
	{\bu}_r({\bf x},t) 
	\equiv \sum_{j=1}^r a_j(t) \bphi_j({\bf x}) \, ,
	\label{eqn:g-rom-1}
\end{equation}
where $\{a_{j}(t)\}_{j=1}^{r}$ are the sought time-varying coefficients, which are determined by solving the following system of equations: $\forall \, i = 1, \ldots, r,$
    \begin{eqnarray}
        \left(
            \frac{\partial \bu_r}{\partial t} , \bphi_{i}
        \right)
        + Re^{-1} \, \left(
            \nabla \bu_r ,
            \nabla \bphi_{i}
        \right)
        + \biggl(
            (\bu_r \cdot \nabla) \, \bu_r ,
            \bphi_{i}
        \biggr)
        = 0 \, .
    \label{eqn:g-rom-weak}
    \end{eqnarray}
In~\eqref{eqn:g-rom-weak}, we assume that the modes $\{ \bphi_1, \ldots, \bphi_r \}$ are perpendicular to the discrete pressure space. 
Plugging~\eqref{eqn:g-rom-1} into~\eqref{eqn:g-rom-weak} yields the {\it Galerkin ROM (G-ROM)}:  
\begin{eqnarray}
	\dot{\ba}
	= A \, \ba
	+ \ba^{\top} \, B \, \ba \, ,
	\label{eqn:g-rom}
\end{eqnarray}
which can be written componentwise as follows: $\forall \, i = 1 \ldots r,$
\begin{equation}
	\dot{a}_i	
	 = \sum_{m=1}^{r} A_{i m} \, a_m(t)
	 + \sum_{m=1}^{r} \sum_{n=1}^{r} B_{i m n} \, a_n(t) \, a_m(t) \, ,
	\label{eqn:g-rom-2}
\end{equation}
where
$A_{im} = - Re^{-1} \, \left( \nabla \bphi_m , \nabla \bphi_i \right)$
and
$B_{imn} = - \bigl( \bphi_m \cdot \nabla \bphi_n , \bphi_i \bigr)$.

The DDF-ROM framework is constructed in two steps.  
In the first step, we propose {\it ROM spatial filtering} to discover the {\it exact mathematical formula} for the ``Correction'' in the Proj-ROM~\eqref{eqn:proj-rom-modified}, which is reminiscent of {\it large eddy simulation (LES)}~\cite{sagaut2006large}.
In this paper, we exclusively use the ROM projection~\cite{wang2012proper,wells2017evolve} as a   spatial filter, but we note that we could also use other spatial filters (e.g., the ROM differential filter~\cite{wells2017evolve,xie2017approximate}).

For a fixed $r \leq d$ and a given $\bu \in \bX^h$ (where $\bX^h$ is the FE space), the ROM projection~\cite{wang2012proper,wells2017evolve} seeks $\obu^r \in \bX^{r}$ such that
        \begin{eqnarray}
            \left( \obu^r, \bphi_j \right)
            = ( \bu , \bphi_j )
            \quad \forall \, j=1, \ldots r \, .
            \label{eqn:rom-projection}
        \end{eqnarray}
Filtering the NSE (see Section 3.2 in~\cite{xie2018data} for details), we obtain the {\it spatially filtered ROM}:
\begin{eqnarray}
	\left(
            \frac{\partial \bu_r}{\partial t} , \bphi_i
        \right)
        + Re^{-1} \, \biggl(
            \nabla \bu_r ,
            \nabla \bphi_i
        \biggr)
        + \biggl(
            \bigl( \bu_r \cdot \nabla \bigr) \, \bu_r ,
            \bphi_i
        \biggr)
        + \biggl(
            \btau_r^{SFS} ,
            \bphi_i
        \biggr)
        = {\bf 0} \, ,
    \label{eqn:f-rom-weak}
\end{eqnarray}
where the ROM stress tensor is
\begin{eqnarray}
\btau_r^{SFS}
= \overline{\bigl( {\bu_d} \cdot \nabla \bigr) \, {\bu_d}}^r 
- \bigl( {\bu_r} \cdot \nabla \bigr) \, {\bu_r} \, .
\label{eqn:les-rom-6}
\end{eqnarray}
The spatially filtered ROM can be written as:  
\begin{eqnarray}
	\dot{\ba}
	= A \, \ba
	+ \ba^{\top} \, B \, \ba 
	+ \btau \, ,
	\label{eqn:f-rom}
\end{eqnarray}
where $A$ and $B$ are the same as those in~\eqref{eqn:g-rom} and the components of $\btau$ are given by 
\begin{eqnarray}
	\tau_i
	= - \biggl(
            \btau_r^{SFS} ,
            \bphi_i
        \biggr) \, ,
        \qquad 
        i = 1, \ldots, r \, .
\label{eqn:les-rom-6b}
\end{eqnarray}

The filtered ROM~\eqref{eqn:f-rom} is an $r$-dimensional ODE system for $\bur$.
Since $r \ll N$, the F-ROM~\eqref{eqn:f-rom} is a computationally efficient surrogate model for the FOM (i.e., the FE approximation of the NSE, which is an $N$-dimensional ODE system).
We emphasize, however, that the F-ROM~\eqref{eqn:f-rom} is not a closed system of equations, since the ROM stress tensor $\btau_r^{SFS}$ (which is used in the definition of $\btau$; see~\eqref{eqn:les-rom-6b}) depends on $\bu_d$ (see~\eqref{eqn:les-rom-6}).
Thus, to close the F-ROM~\eqref{eqn:f-rom}, we need to solve the {\it ROM closure problem}~\cite{chekroun2017markovian,galletti2007accurate,lu2017data,noack2008finite,sirisup2004spectral,wang2012proper}, i.e., to find a formula $ \btau \approx \btau(\ba)$.  
The {\it explicit formula for $\btau$} (see~\eqref{eqn:les-rom-6} and \eqref{eqn:les-rom-6b}), {\it  allows for the first time the use of data-driven modeling of the entire missing ROM information}. \\[-0.3cm]

In the second step of the DDF-ROM construction, we use data-driven modeling to solve the ROM closure problem, i.e., to find a formula $ \btau \approx \btau(\ba)$ in~\eqref{eqn:f-rom}.
To make the F-ROM~\eqref{eqn:f-rom} resemble the standard G-ROM~\eqref{eqn:g-rom}, we make the following ansatz: 
\begin{eqnarray}
     \btau(\ba)
     	\approx \btau^{ansatz}(\ba)
	= \tA \, \ba
	+ \ba^{\top} \, \tB \, \ba 
 \, . 
	\label{eqn:c-rom-6}
\end{eqnarray}
Using ansatz~\eqref{eqn:c-rom-6} in the F-ROM~\eqref{eqn:f-rom} yields a closed system of equations.
To find $\tA$ and $\tB$ in~\eqref{eqn:c-rom-6}, we use data-driven modeling.
That is, we find $\tA$ and $\tB$ that ensure the highest accuracy of the vector $\btau$ in the filtered ROM~\eqref{eqn:f-rom}.
To this end, we minimize the $L^2$-norm of the difference between $\btau$ computed with the FOM data and equations~\eqref{eqn:les-rom-6} and~\eqref{eqn:les-rom-6b}, and $\btau$ computed with the ansatz~\eqref{eqn:c-rom-6} and the ROM coefficients obtained from the snapshots.
Thus, we solve the following optimization problem~\cite{noack2005need,peherstorfer2016data}:
\begin{eqnarray}
	\min_{\substack{\tA \in \R^{r \times r} \\[0.1cm] \tB \in \R^{r \times r \times r}}} \, 
	\sum_{j = 1}^{M} \| \btau^{true}(t_j) - \btau^{ansatz}(t_j) \|^2
	\label{eqn:ddf-rom-least-squares} \, ,
\end{eqnarray}
where $\| \cdot \|$ is the Euclidian norm in $\R^r$ and $\btau^{true}(t_j)$ is the {\it true} $\btau(t_j)$ computed from the snapshot data.
The optimal $\tA^{opt}$ and $\tB^{opt}$ are used in the filtered ROM~\eqref{eqn:f-rom}, yielding the {\it data-driven filtered ROM (DDF-ROM)}
\begin{equation}
	\boxed{
	\dot{\ba}
	= \left( A + \tA \right) \, \ba
	+ \ba^{\top} \, \left( B + \tB \right) \, \ba \, . 
	}
	\label{eqn:ddf-rom}
\end{equation}

\begin{figure}[h!]
	\centering
	\includegraphics[width = 1\textwidth, height=.3\textwidth, clip]{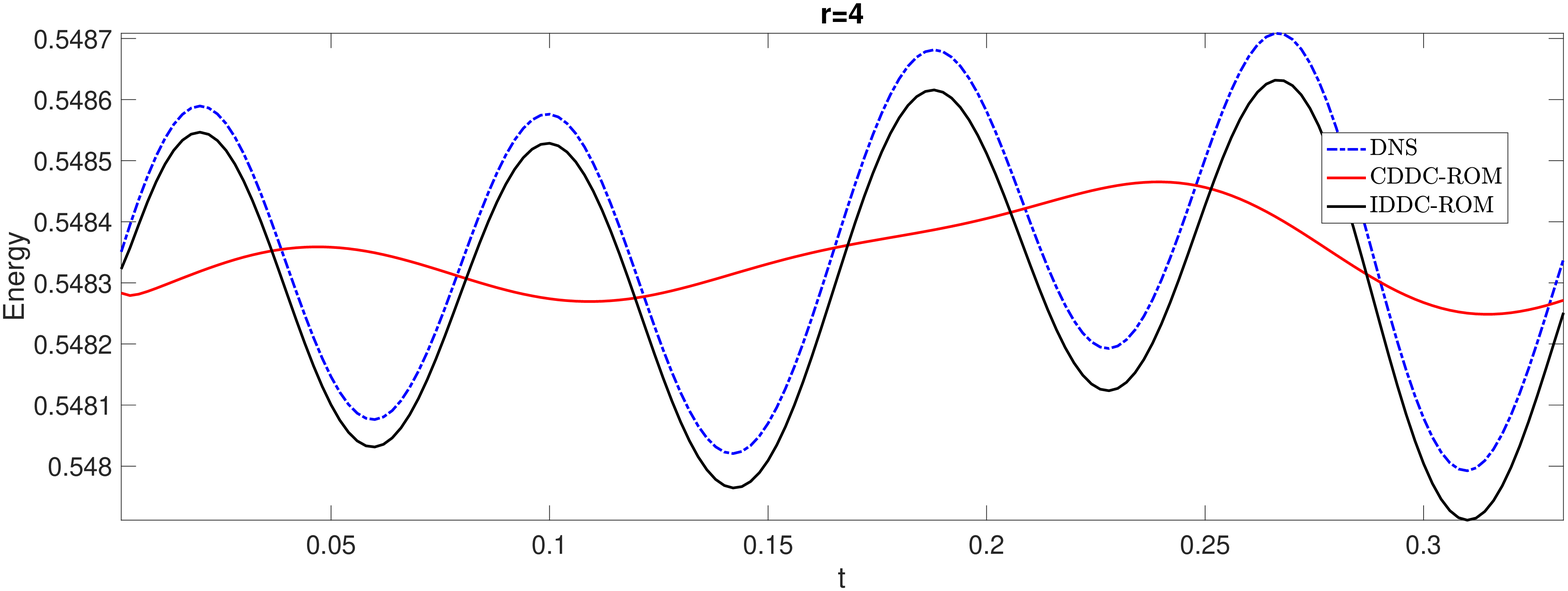}
	\caption{
	Plots of energy coefficients vs. time for  the ``ideal'' DDF-ROM, DNS, and CDDF-ROM for flow past a circular cylinder with $Re=100$ (see Section~\ref{sec:numerical-results} for details).
	\label{fig:iddf-rom}
	}
\end{figure}

\begin{remark}[The Role of the Ansatz in the DDF-ROM: \ The Ideal DDF-ROM]
The effect of the ansatz~\eqref{eqn:c-rom-6} on the numerical results (see Section~\ref{sec:numerical-results} for details) is illustrated in Fig.~\ref{fig:iddf-rom}. 
An ``ideal'' DDF-ROM (i.e., a DDF-ROM that uses the true $\btau$) is almost as accurate as the FOM, even with $r=4$!
Using the ansatz (and solving the least squares problem) yields results that are qualitatively similar, but loses some of the ``ideal'' DDF-ROM's accuracy.
(Of course, we note that the ``ideal'' DDF-ROM is used only for illustration purposes; it is not a  practical ROM, since it can only be used with the training data.)
\end{remark}


\section{Physically Constrained DDF-ROM (CDDF-ROM)}
	\label{sec:cddf-rom}

The operators $A$ and $B$ in the G-ROM~\eqref{eqn:g-rom} satisfy several physical constraints:
First, the operator $A$ is negative semidefinite: 
\begin{eqnarray}
	\ba^{\top} A \, \ba
	=
	- \| \nabla \bur \|^{2}
	\leq 
	0 \, .
	\label{eqn:a-negative-semidefinite}
\end{eqnarray}
Second, if one uses the skew-symmetric formulation of the nonlinearity, the operator $B$ satisfies 
\begin{eqnarray}
	\ba^{\top} \left[ \ba^{\top} \, B \, \ba \right] 
	=
	0 .
	\label{eqn:b-conservative}
\end{eqnarray}
A natural question is whether the DDF-ROM operators $\tA$ and $\tB$ should also satisfy any physical constraints.

We {\it conjecture} that the DDF-ROM operators $\tA$ and $\tB$ satisfy physical constraints that are consistent with the physical constraints satisfied by the ROM subfilter-scale stress tensor $\btau$, which is used in their construction (see equations \eqref{eqn:les-rom-6} and \eqref{eqn:les-rom-6b}).
In~\cite{CSB03}, it was shown that the role of the ROM subfilter-scale stress tensor (i.e., the ``Correction'' in~\eqref{eqn:proj-rom-modified}) is to dissipate energy.
Thus, we conjecture that the constraint that we need to enforce in the DDF-ROM construction is
\begin{eqnarray}
			\ba^{\top} [ \tA \, \ba + \ba^{\top} \, \tB \, \ba ]
			\leq
			0 \, .
		\label{eqn:ddf-rom-constraint-atilde-btilde}
	\end{eqnarray}
The challenge in enforcing this constraint is that $\tA$ yields quadratic terms, whereas $\tB$ yields cubic terms.
Thus, we propose to replace~\eqref{eqn:ddf-rom-constraint-atilde-btilde} with the following constraints:
	\begin{eqnarray}
		\boxed{
			\ba^{\top} \tA \, \ba
			\leq 
			0 
			\quad \text{and} \quad  
			\ba^{\top} \left[ \ba^{\top} \, \tB \, \ba \right]
			=
			0 \, ,
		}
		\label{eqn:ddf-rom-constraint-atilde-btilde-separate}
	\end{eqnarray}
which are easier to implement.
Furthermore, the constraints~\eqref{eqn:ddf-rom-constraint-atilde-btilde-separate} resemble the constraints satisfied by the operators $A$ and $B$ in the G-ROM, i.e., the constraints~\eqref{eqn:a-negative-semidefinite} and~\eqref{eqn:b-conservative}, respectively.
To enforce the first constraint in~\eqref{eqn:ddf-rom-constraint-atilde-btilde-separate}, sufficient conditions are~\cite{kondrashov2015data}:
\begin{eqnarray}		
	\widetilde{A}_{i i}
	\leq
	0 ,
	\ 
	\forall \, i = 1, \ldots, r, 
	\label{eqn:atilde-constraint-componentwise-1}
\end{eqnarray}
and 
\begin{eqnarray}
	\widetilde{A}_{i j}
	= 
	- \widetilde{A}_{j i} ,
	\ 
	\forall \, i,j = 1, \ldots, r, \ i \neq j \, .
	\label{eqn:atilde-constraint-componentwise-2}
\end{eqnarray}
To enforce the second constraint in~\eqref{eqn:ddf-rom-constraint-atilde-btilde-separate}, sufficient conditions are:
\begin{eqnarray}
	\widetilde{B}_{i i i}
	&=&
	0 , 
	\quad
	\forall \, i = 1, \ldots, r,
	\label{eqn:btilde-constraint-componentwise-1} \\[0.3cm]
	\widetilde{B}_{i i j}
	+ 
	\widetilde{B}_{i j i}
	+
	\widetilde{B}_{j i i}
	&=&
	0 , 
	\quad 
	\forall \, i,j = 1, \ldots, r, \ i \neq j \, ,
	\label{eqn:btilde-constraint-componentwise-2} \\[0.3cm]
	\hspace*{-0.7cm}
	\widetilde{B}_{i j k}
	+ 
	\widetilde{B}_{i k j}
	+
	\widetilde{B}_{j i k}
	+
	\widetilde{B}_{j k i}
	+
	\widetilde{B}_{k i j}
	+
	\widetilde{B}_{k j i}
	&=&
	0 , 
	\quad 
	\forall \, i,j,k = 1, \ldots, r, \ i \neq j \neq k \neq i \, .
	\label{eqn:btilde-constraint-componentwise-3}
\end{eqnarray}

We emphasize that, from the implementation standpoint, enforcing constraints~\eqref{eqn:ddf-rom-constraint-atilde-btilde-separate} in the data-driven ROM closure modeling is straightforward:
Instead of solving the unconstrained least squares problem~\eqref{eqn:ddf-rom-least-squares}, we just need to solve a {\it constrained least squares problem} with constraints~\eqref{eqn:ddf-rom-constraint-atilde-btilde-separate}. 
%
%
Thus, the {\it physically-constrained data-driven filtered ROM (CDDF-ROM)} has the form
\begin{equation}
	\boxed{
	\dot{\ba}
	= \left( A + \tA \right) \, \ba
	+ \ba^{\top} \, \left( B + \tB \right) \, \ba \, ,
	}
	\label{eqn:cddf-rom}
\end{equation}
where the operators $\tA$ and $\tB$, instead of solving the unconstrained least squares problem~\eqref{eqn:ddf-rom-least-squares} (as is done in the DDF-ROM), solve the following {\it constrained} least squares problem:
\begin{eqnarray}
	\boxed{
	\min_{
		\substack{\tA \in \R^{r \times r} \\[0.1cm] 
				\tB \in \R^{r \times r \times r} \\[0.1cm] 
				\ba^{\top} \tA \, \ba \leq 0 \\[0.1cm] 
				\ba^{\top} \left[ \ba^{\top} \, \tB \, \ba \right] = 0
				}
		} \, 
	\sum_{j = 1}^{M} \| \btau^{true}(t_j) - \btau^{ansatz}(t_j) \|^2
	\label{eqn:cddf-rom-least-squares} \, .
	}
\end{eqnarray}

\bigskip

We note that, in data-driven modeling, physical constraints have been enforced in, e.g.,~\cite{kondrashov2015data,loiseau2018constrained}.
We emphasize, however, that the CDDF-ROM setting is different from that used in~\cite{kondrashov2015data,loiseau2018constrained}:
Indeed, in~\cite{kondrashov2015data,loiseau2018constrained} the authors consider {\it pure} data-driven ROMs~\eqref{eqn:dd-rom}, which only involve operators $\tA$ and $\tB$ (without $A$ and $B$).
The new DDF-ROM~\eqref{eqn:ddf-rom}, on the other hand, is a {\it hybrid} data-driven/projection ROM, so it involves both G-ROM operators ($A$ and $B$) and data-driven operators ($\tA$ and $\tB$).


\section{Numerical Results}
	\label{sec:numerical-results}

In this section, we perform a numerical investigation of the new CDDF-ROM~\eqref{eqn:cddf-rom}--\eqref{eqn:cddf-rom-least-squares}.
Specifically, we investigate whether the constrained data-driven modeling in the new CDDF-ROM construction yields more accurate results than the unconstrained data-driven modeling in the DDF-ROM proposed in~\cite{xie2018data}.
To this end, we compare the new CDDF-ROM with the DFF-ROM in the numerical simulation of a 2D flow past a circular cylinder at Reynolds numbers $Re=100, Re=500$, and $Re=1000$~\cite{ST96,mohebujjaman2017energy}.
As a benchmark for our comparison, we use the FOM data obtained with a FE simulation.
Since most of the numerical investigation in this section is performed for $Re=100$, we present the background material (e.g., computational setting and ROM construction) for this case and only list the main differences from the $Re=100$ test case when we consider the higher Reynolds number cases (i.e., $Re=500$ and $Re=1000$). 
The rest of the section is organized as follows:
In Section~\ref{sec:test-problem-setup}, we describe the test problem setup.
In Section~\ref{sec:snapshot-rom-generation}, we outline the snapshot and ROM generation.
In Section~\ref{sec:computational-efficiency}, we discuss the computational efficiency of the new CDDF-ROM.
In Section~\ref{sec:cddf-rom-vs-ddf-rom}, we compare the new CDDF-ROM with the DDF-ROM.
Finally, in Section~\ref{sec:cross-validation}, we present a {\it cross-validation} of the two ROMs, i.e., we test the ROMs for data that was {\it not} used to train the ROM closure model.
Specifically, we investigate the {\it predictive} capabilities of the CDDF-ROM.

\subsection{Test Problem Setup}
	\label{sec:test-problem-setup}
	
\begin{figure}[h!]
\begin{center}
\includegraphics[width=0.7\textwidth,height=0.24\textwidth, trim=0 0 0 0, clip]{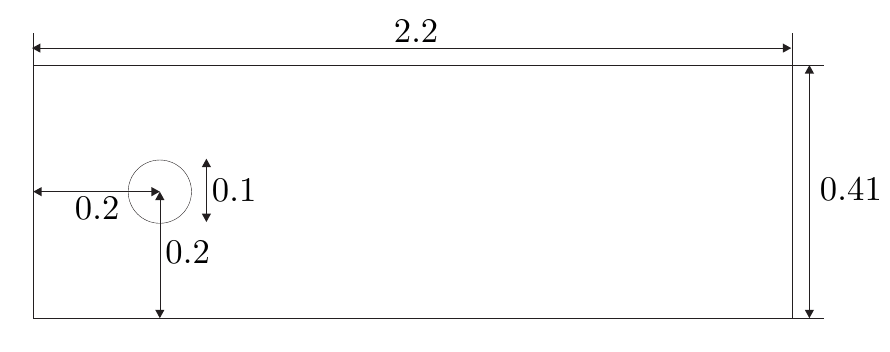}
\end{center}
\caption{\label{cyldomain} Channel flow around a cylinder domain.}
\end{figure}

The domain is a $2.2\times 0.41$ rectangular channel with a radius=$0.05$ cylinder, centered at $(0.2,0.2)$, see Figure~\ref{cyldomain}.  
No slip boundary conditions are prescribed for the walls and on the cylinder, the inflow profile is given by~\cite{john2004reference,rebholz2017improved} $u_{1}(0,y,t)=u_{1}(2.2,y,t)=\frac{6}{0.41^{2}}y(0.41-y) \, , u_{2}(0,y,t)=u_{2}(2.2,y,t)=0$, and appropriate outflow boundary conditions are used.
The kinematic viscosity is $\nu=10^{-3}$, there is no forcing, and the flow starts from rest.

\subsection{Snapshot and ROM Generation}
	\label{sec:snapshot-rom-generation}


To compute the snapshots, we use the commonly used linearized BDF2 temporal discretization,  together with a FE spatial discretization utilizing the Scott-Vogelius element.  
On time step $1$, we use a backward Euler temporal discretization.  
All simulations use a time step size of $\Delta t=0.002$.
\blue{
We compute until a statistically steady state is reached (which occurs at about $T=5$), and then compute to $T=17$.
}
Snapshots are taken to be the solutions at each time step from $T=7$ to $T=7.332$, which corresponds to one period. Thus, in total 166 snapshots were collected.
We compute on two different meshes, which provide approximately 103K and 35K velocity degrees of freedom.  
The 103K mesh gives essentially a fully resolved solution, and the lift and drag predictions agree well with results from fine discretizations in~\cite{caiazzo2014numerical,ST96}:
$c_{d,max} = 3.2261,\ c_{l,max} = 1.0040.$
Results from the 35K meshes are only slightly less accurate.
The results for the $Re=100$ test case are obtained on the 35K mesh; the results for the $Re=500$ and $Re=1000$ test cases are obtained on the 103K mesh. 

\begin{figure}[h!]
	 	 	\centering
	 	 	\includegraphics[width = .9\textwidth, height=.3\textwidth,viewport=0 25 1000 430, clip]{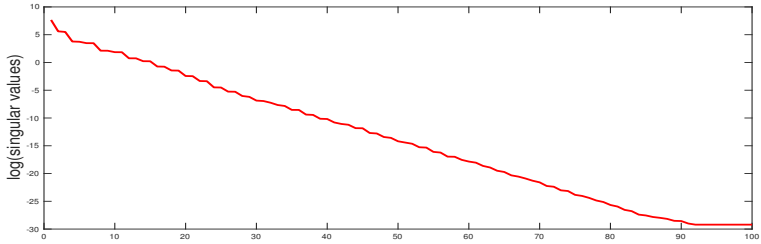}
	 	 	\caption{
			Plots of singular values  vs. index, for flow past a cylinder with $Re=100$.
			\label{fig:singular-values} 
			}
\end{figure}

The ROM modes are created from the snapshots in the usual way.  
The first mode is chosen to be the snapshot average, which satisfies the boundary conditions.  
This mode is then subtracted from the snapshots, and finally an eigenvalue problem is solved to find the dominant modes of these adjusted snapshots (see \cite{caiazzo2014numerical} for a more detailed description of the process). 
The singular values of the snapshot matrix are plotted in Figure~\ref{fig:singular-values}.

With the dominant modes created, the ROM is constructed as discussed in Section~\ref{sec:ddf-rom} using the BDF2 temporal discretization.  
In all of our tests, just as in the FE simulations, we take $\Delta t=0.002$; this choice creates no significant temporal error in any of our simulations (tests were done with varying $\Delta t$ to verify that 0.002 is sufficiently small).  
The ROM initial condition at $T=6.998$ is the $L^2$ projection of the FE solution at $T=6.998$ into the ROM space.
The ROM initial condition at $T=7$ is obtained by using the backward Euler method.
The ROMs are run from this initial time (now called $t=0$), and continued to $t=10$.
The ROMs are tested using three $r$ values: $r=3, r=4$, and $r=6$.

\vspace*{0.3cm}

To allow more flexibility in the computational implementation, instead of using the constraints $\tilde{A}_{ii} \leq0$ in~\eqref{eqn:atilde-constraint-componentwise-1}, we introduce a parameter $\epsilon \geq 0$ (which is typically small) and enforce the constraints $\tilde{A}_{ii} \leq -\epsilon$. 
To solve the constrained optimization problem, we used MATLAB toolbox lsqlin with interior-point method, where $\texttt{ConstrainedTolerence}=10^{-9}$, $\texttt{OptimalityTolerance}=10^{-9}$, and $\texttt{StepTolerance}$ $=10^{-9}$. 

\subsection{Computational Efficiency}
	\label{sec:computational-efficiency}

Although the CDDF-ROM and DDF-ROM are more accurate than the G-ROM, the computational cost of calculating $\tA$ and $\tB$ in the offline phase can be significant.
Thus, in~\cite{xie2018data}, we proposed the following practical approach for reducing the computational cost of the $\tA$ and $\tB$ calculation:
Since $d$, the rank of the snapshot matrix, can be large in practical applications,instead of using $\bu_d \in \bX^d$ to compute the Correction term in~\eqref{eqn:les-rom-6}, we utilized the following approximation:
\begin{eqnarray}
	- \left(
        		\overline{\bigl( {\bu_d} \cdot \nabla \bigr) \, {\bu_d}}^r 
		- \bigl( {\bu_r} \cdot \nabla \bigr) \, {\bu_r}  ,
            \bphi_i
        \right) 
        \approx 
        - \left(
        		\overline{\bigl( {\bu_m} \cdot \nabla \bigr) \, {\bu_m}}^r 
		- \bigl( {\bu_r} \cdot \nabla \bigr) \, {\bu_r}  ,
            \bphi_i
        \right) ,
        \qquad
        \label{eqn:ddf-rom-efficiency}
\end{eqnarray}
$\forall \, i = 1, \ldots, r$.
In~\eqref{eqn:ddf-rom-efficiency}, we replaced $\bu_d$ with $\bu_{m}$ and the ROM projection on $\bX^d$ with the ROM projection on $\bX^{m}$.
In practice, the parameter $m$ in~\eqref{eqn:ddf-rom-efficiency}, where $r \leq m \leq d$, should be chosen to balance accuracy and efficiency in the DDF-ROM~\cite{xie2018data}.
In our numerical investigation, however, we varied the parameter value $m$ to achieve the highest level of accuracy.

\subsection{Ill-Conditioning}
	\label{sec:ill-conditioning}

In Section 4 in~\cite{xie2018data}, we observed that the least squares problem used to compute the DDF-ROM operators $\tA$ and $\tB$ was ill-conditioned.
We also noted that, in general, data-driven least squares problems can be ill-conditioned; e.g., the least squares problem in the data-driven operator inference method proposed in~\cite{peherstorfer2016data} was also ill-conditioned.
To remedy this ill-conditioning, in Algorithm 1 in~\cite{xie2018data}, we used the truncated singular value decomposition (SVD)~\cite{demmel1997applied} (see~\cite{peherstorfer2016data} for alternative approaches).
The CDDF-ROM's constrained least squares problem~\eqref{eqn:cddf-rom-least-squares} is again ill-conditioned, just as the DDF-ROM's unconstrained least squares problem.
To tackle this ill-conditioning, we use again the truncated SVD procedure (see Step 6 of Algorithm 1 in~\cite{xie2018data}).

\subsection{CDDF-ROM vs DDF-ROM}
	\label{sec:cddf-rom-vs-ddf-rom}

In this section, we compare the new CDDF-ROM with the DDF-ROM.
As a benchmark, we use the FOM.
We present results for three Reynolds numbers: $Re=100, Re=500$, and $Re=1000$.

Just as the DDF-ROM, the CDDF-ROM is sensitive to variations in its parameters: 
$tol$ (the tolerance used in the truncated-svd algorithm) and $d$ (the dimension of the projection space).
In addition, CDDF-ROM is also sensitive to $\epsilon$, which is the parameter used in the constrained least squares problem in Section~\ref{sec:cddf-rom}.
In what follows, we present results for the optimal parameter values that were found by trial and error.


\subsubsection{Reynolds Number ${\bf Re=100}$}
In this section, we present results for $Re=100$.

\begin{figure}[h]
  \includegraphics[width=.9\textwidth,viewport=0 25 1350 460, clip]{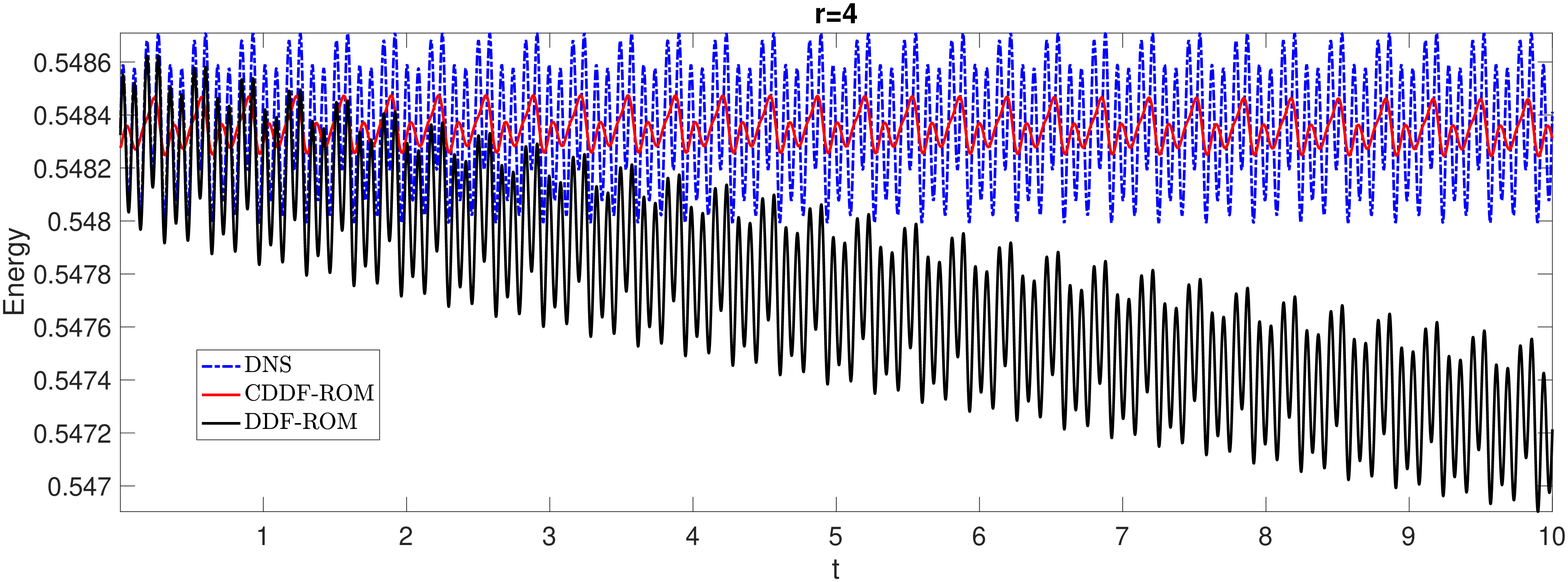}
  \includegraphics[width=.9\textwidth,viewport=0 25 1350 437, clip]{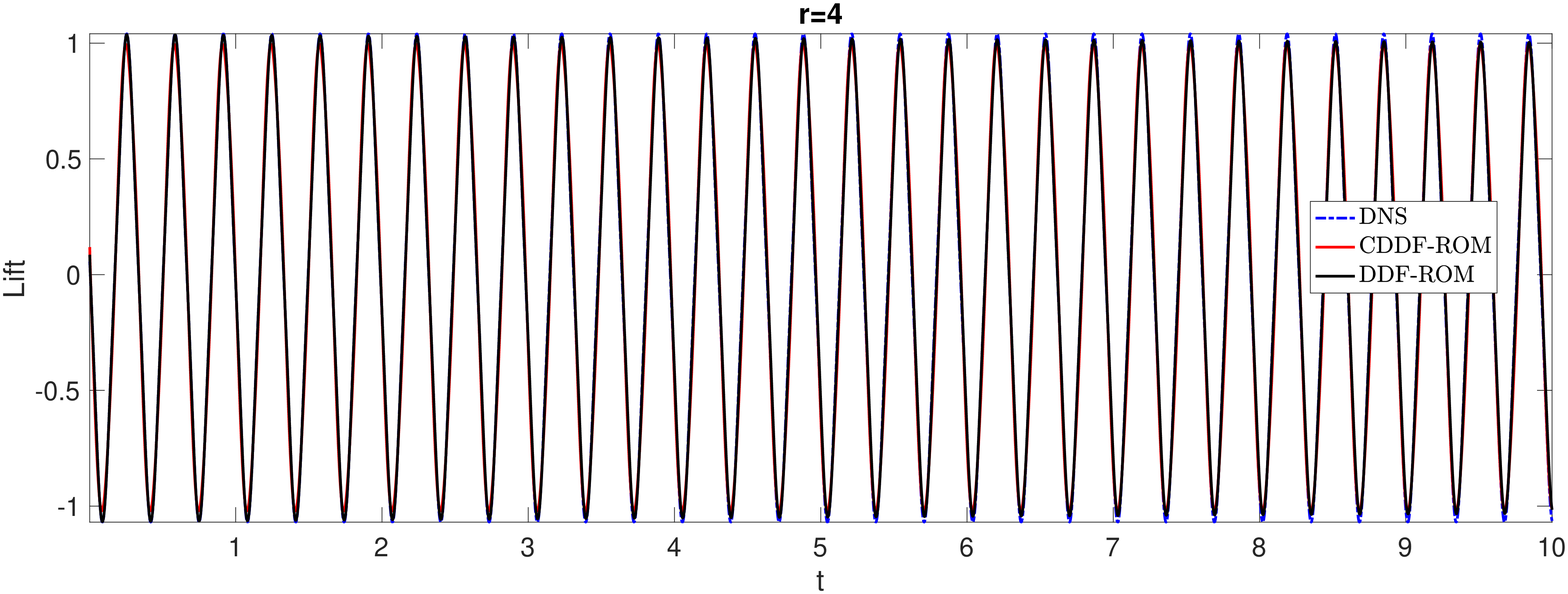}\\
  \includegraphics[width=.9\textwidth,viewport=0 25 1350 437, clip]{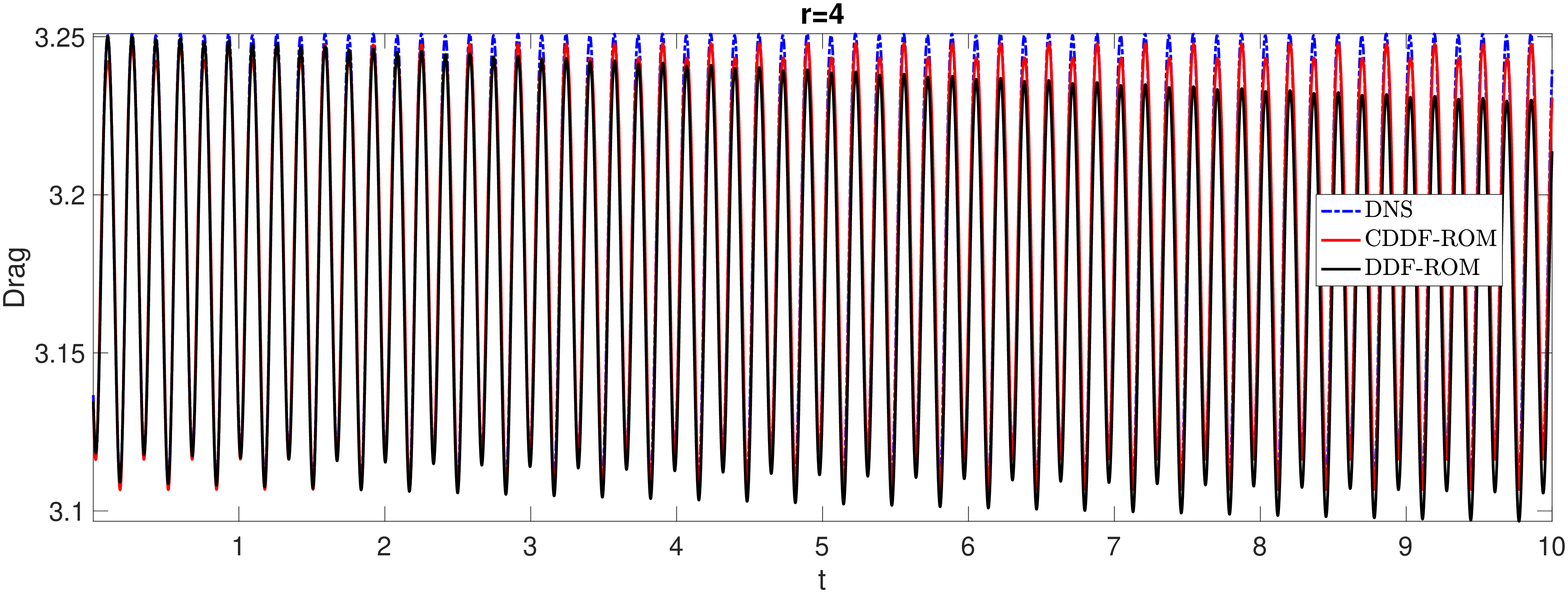}
     \caption{
     	Flow past a circular cylinder with $Re=100$:
     	Plots of energy, lift, and drag coefficients vs. time for DNS, DDF-ROM ($r=4$), and CDDF-ROM ($r=4$).
	\label{cddfr4d7}
	}
 \end{figure}

First, we consider $r=4$.
We note that $m=r+3$ is the minimum $m$ value for which CDDF-ROM yields accurate results. 
Lower $d$ values yield inaccurate results in our numerical investigation.
For the CDDF-ROM, we use $tol=7\times 10^{-3}, \epsilon=7.1\times 10^{-10}$; for the DDF-ROM, we use $tol=7\times 10^{-3}$.
The energy evolution in Fig.~\ref{cddfr4d7} clearly shows that the CDDF-ROM is dramatically more accurate than the DDF-ROM.
The same holds for the drag evolution, although the difference is not as big.
Finally, the lift evolution of CDDF-ROM is similar to that of DDF-ROM.
 
\begin{figure}[h]
     \includegraphics[width=.9\textwidth,viewport=0 25 1350 460, clip]{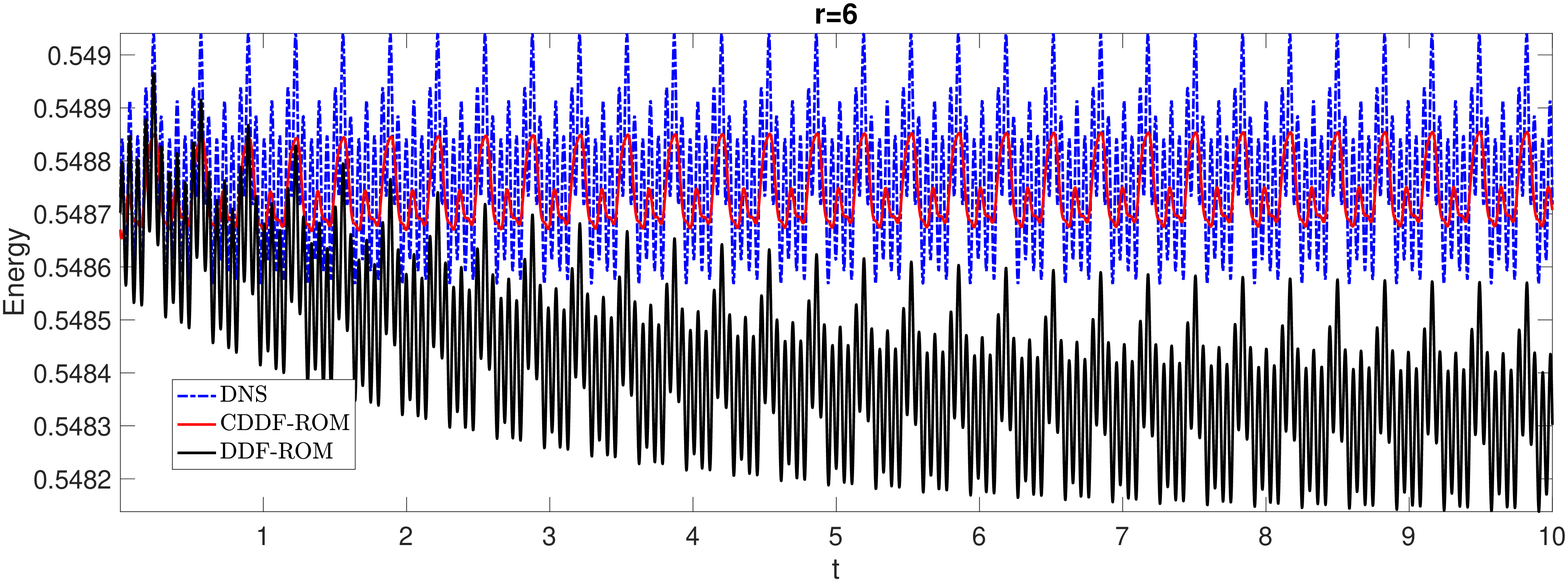}
     \includegraphics[width=.9\textwidth,viewport=0 25 1350 430, clip]{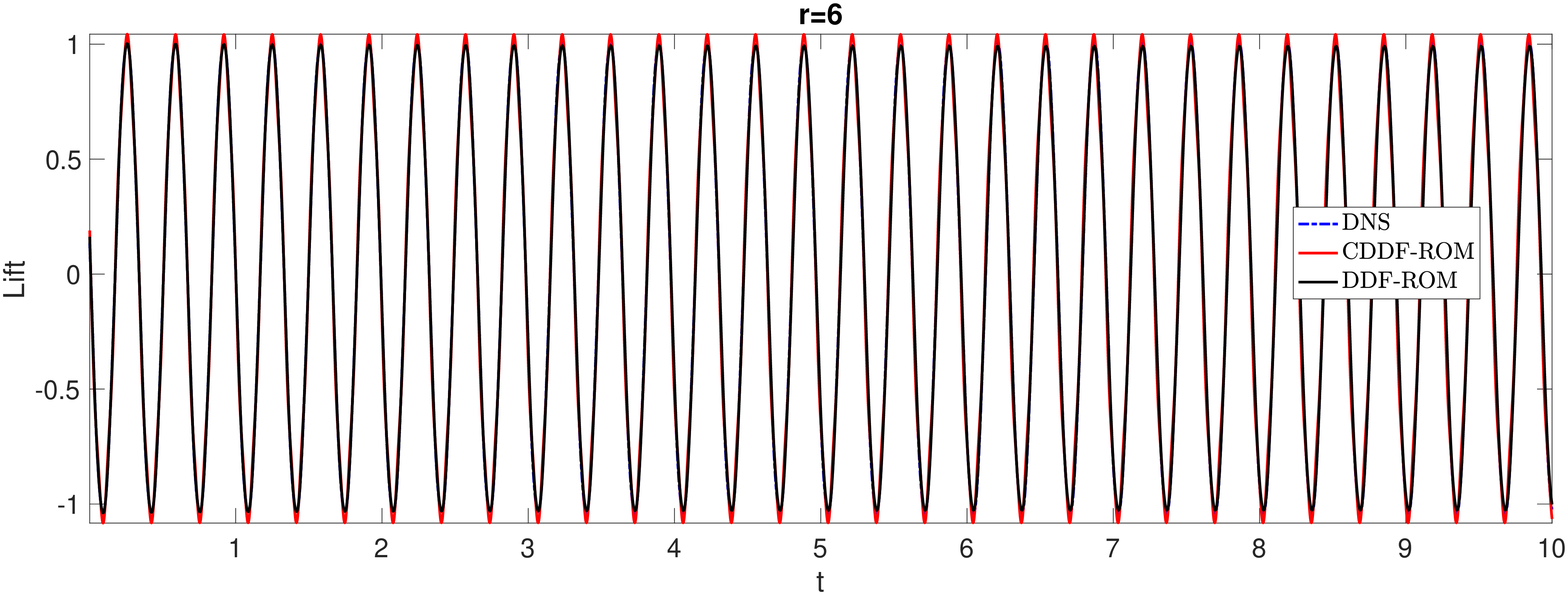}\\
     \includegraphics[width=.9\textwidth,viewport=0 25 1350 435, clip]{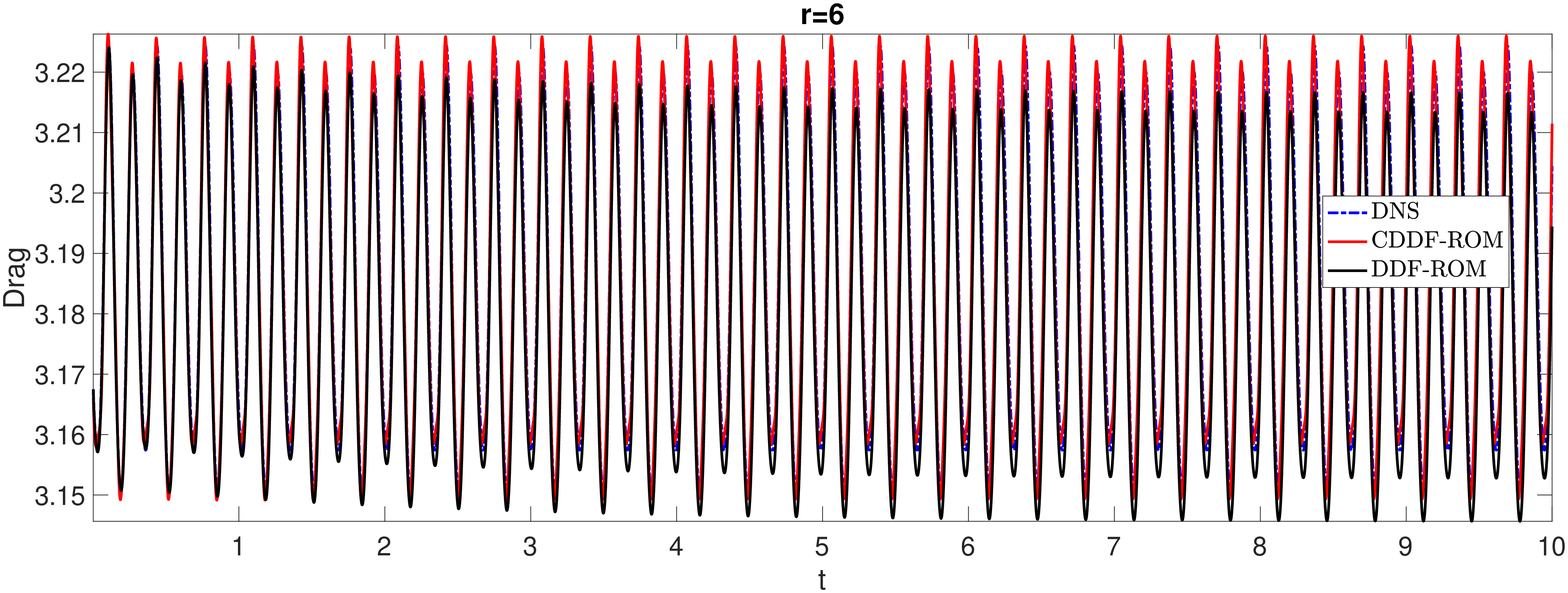}
     \caption{
     	Flow past a circular cylinder with $Re=100$:
     	Plots of energy, lift, and drag coefficients vs. time for DNS, DDF-ROM ($r=6$), and CDDF-ROM ($r=6$).
	\label{dd_vs_cddfr6d7}
	}
    \end{figure}

Next, we investigate the case $r=6$.
This time, the minimum $m$ value that yields accurate results for CDDF-ROM and DDF-ROM is $m = r+1$.
For the CDDF-ROM, we use $tol=1.2\times 10^{-2}, \epsilon=8.5\times 10^{-3}$; for the DDF-ROM, we use $tol=7\times 10^{-3}$.
In Fig.~\ref{dd_vs_cddfr6d7}, we plot the evolution of energy, lift, and drag coefficients for CDDF-ROM, DDF-ROM, and FOM.
As in the $r=4$ case, the CDDF-ROM is significantly more accurate than the DDF-ROM.

We conclude that the CDDF-ROM is significantly more accurate than the DDF-ROM for low $r$ values, for which the latter does not perform well.
We note, however, for higher $r$ values for which the DDF-ROM performs well, the CDDF-ROM does not show a visible improvement.

\subsubsection{Reynolds Number ${\bf Re=500}$}
In this section, we present results for $Re=500$.

For the $Re=500$ test case, to compute the snapshots, we use the same approach as that described in Section~\ref{sec:snapshot-rom-generation}, except that the snapshots are the solutions at each time step from $T=10$ to $T=10.438$, which correspond to one period.
Thus, in this case we collected a total of $219$ snapshots. 
As ROM initial condition, we use the $L^2$ projection of the snapshot at $T=10$ on the ROM space.

\begin{figure}[h]
	\includegraphics[width=.9\textwidth,viewport=0 25 1350 455, clip]{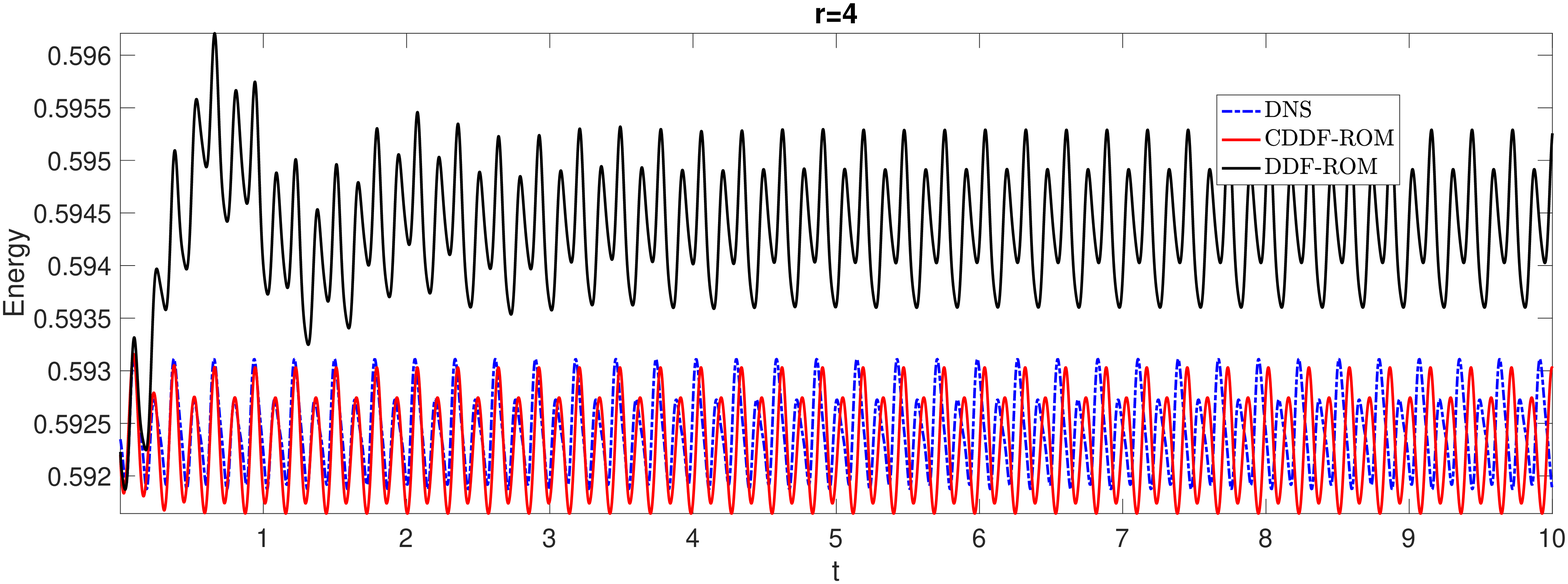}
     \caption{
     	Flow past a circular cylinder with $Re=500$:
     	Plots of energy, lift, and drag coefficients vs. time for DNS, DDF-ROM ($r=4$), and CDDF-ROM ($r=4$).
	\label{Re500r4}
	}
\end{figure}

First, we consider $r=4$; $m=r+1$ is the minimum $m$ value for which CDDF-ROM yields accurate results. 
For the CDDF-ROM, we use $tol=7\times 10^{-3}, \epsilon=0$; for the DDF-ROM, we use $tol=7.5\times 10^{-3}$.
The energy evolution in Fig.~\ref{Re500r4} clearly shows that the CDDF-ROM is dramatically more accurate than the DDF-ROM.
Since the lift and drag evolutions are similar to those for $Re=100$, in what follows we will not include plots for these quantities.

\begin{figure}[h]
	\includegraphics[width=.9\textwidth,viewport=0 25 1350 455, clip]{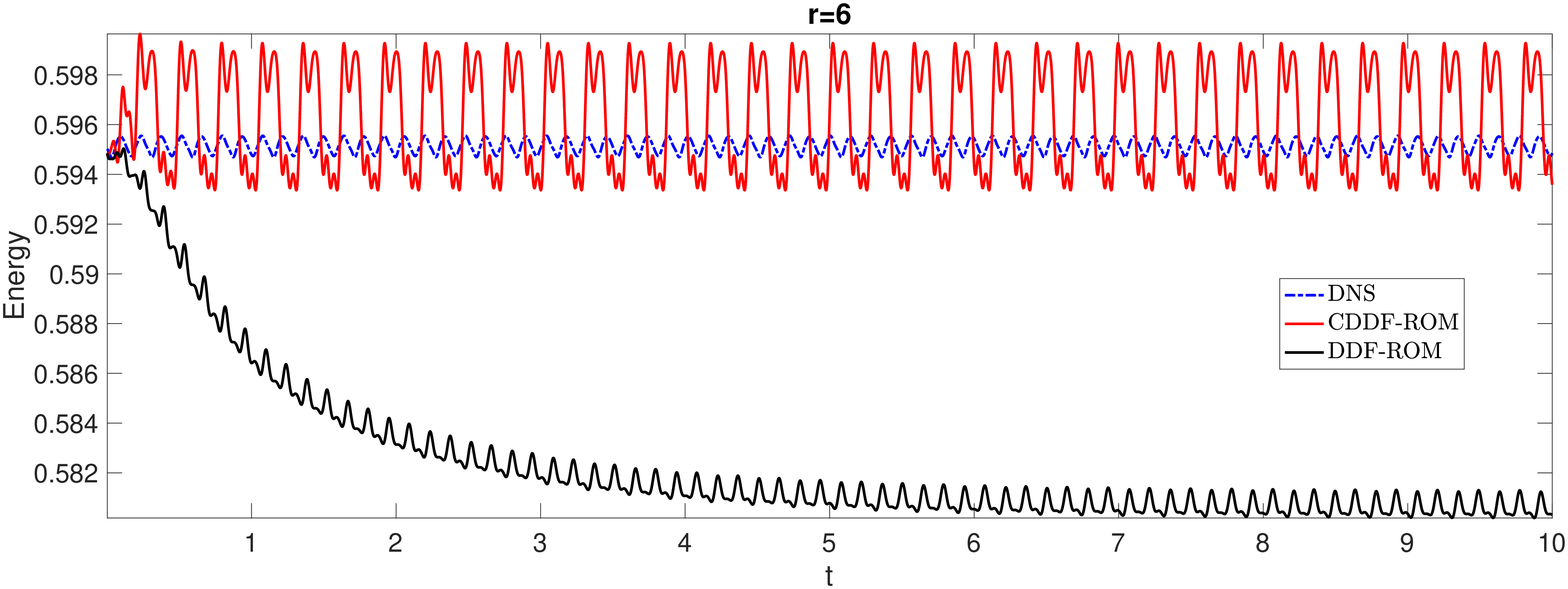}
     \caption{
     	Flow past a circular cylinder with $Re=500$:
     	Plots of energy, lift, and drag coefficients vs. time for DNS, DDF-ROM ($r=6$), and CDDF-ROM ($r=6$).
	\label{Re500r6}
	}
\end{figure}

Next, we investigate the case $r=6$.
Again, the minimum $m$ value that yields accurate results for CDDF-ROM and DDF-ROM is $m = r+1$.
For the CDDF-ROM, we use $tol= 10^{-1}, \epsilon=4\times 10^{-3}$; for the DDF-ROM, we use $tol=5\times 10^{-3}$.
In Fig.~\ref{Re500r6}, we plot the evolution of energy coefficients for CDDF-ROM, DDF-ROM, and FOM.
As in the $r=4$ case, the CDDF-ROM is significantly more accurate than the DDF-ROM.

\subsubsection{Reynolds Number ${\bf Re=1000}$}
In this section, we present results for $Re=1000$.

For the $Re=1000$ test case, to compute the snapshots, we use the same approach as that described in Section~\ref{sec:snapshot-rom-generation}, except that the snapshots are the solutions at each time step from $T=5$ to $T=5.268$, which correspond to one period.
Thus, in this case we collected a total of $134$ snapshots. 
As ROM initial condition, we use the $L^2$ projection of the snapshot at $T=T=5.004$ on the ROM space.

\begin{figure}[h]
	\includegraphics[width=.9\textwidth,viewport=0 25 1350 455, clip]{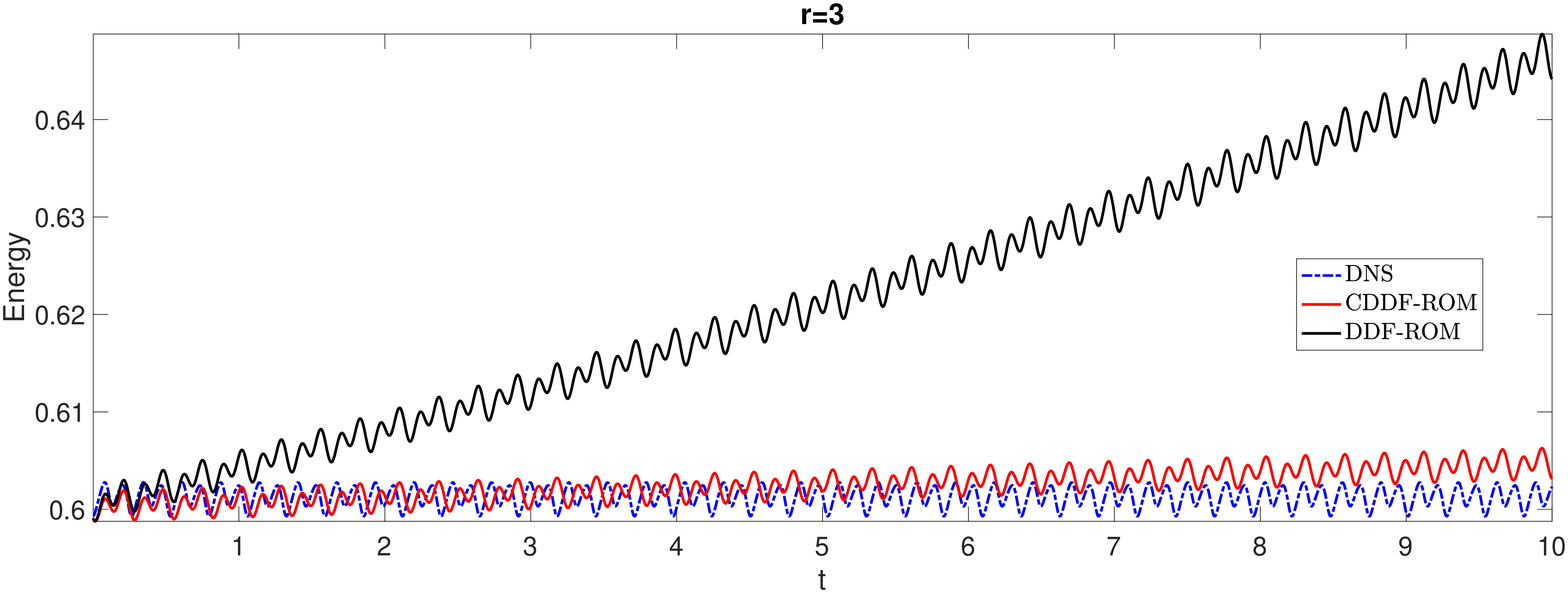}\\
	\includegraphics[width=.9\textwidth,viewport=0 25 1350 437, clip]{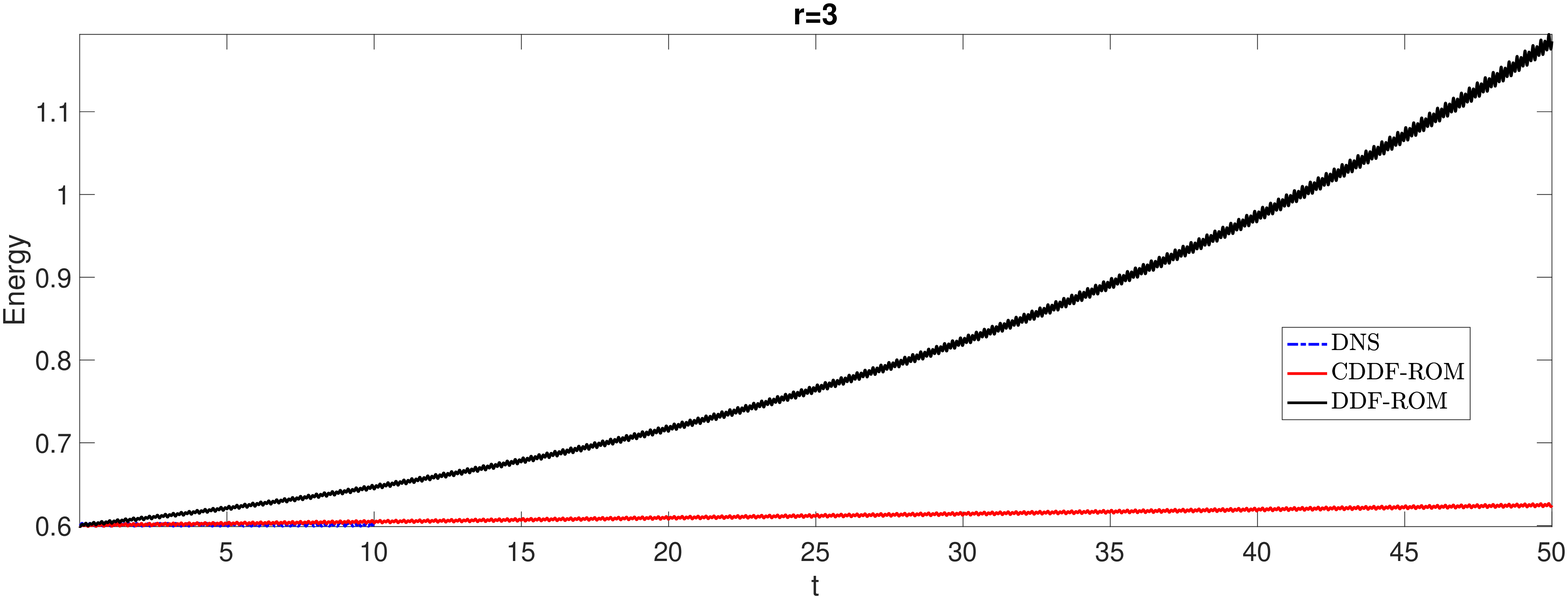}
     \caption{
     	Flow past a circular cylinder with $Re=1000$:
     	Plots of energy coefficients vs. time for DNS, DDF-ROM ($r=3$), and CDDF-ROM ($r=3$).
	\label{fig:Re1000r4}
	}
\end{figure}

First, we consider $r=4$; $m=r+1$ is the minimum $m$ value for which CDDF-ROM yields accurate results. 
For the CDDF-ROM, we use $tol= 0.42, \epsilon=2.5\times 10^{-3}$; for the DDF-ROM, we use $tol= 5\times 10^{-3}$.
The energy evolution in Fig.~\ref{fig:Re1000r4} clearly shows that the CDDF-ROM is dramatically more accurate than the DDF-ROM.
The plot at the bottom of Fig.~\ref{fig:Re1000r4} also shows that although the CDDF-ROM energy grows in time, it does so at a slower rate than the DDF-ROM energy.

\begin{figure}[h]
	\includegraphics[width=.9\textwidth,viewport=0 25 1350 452, clip]{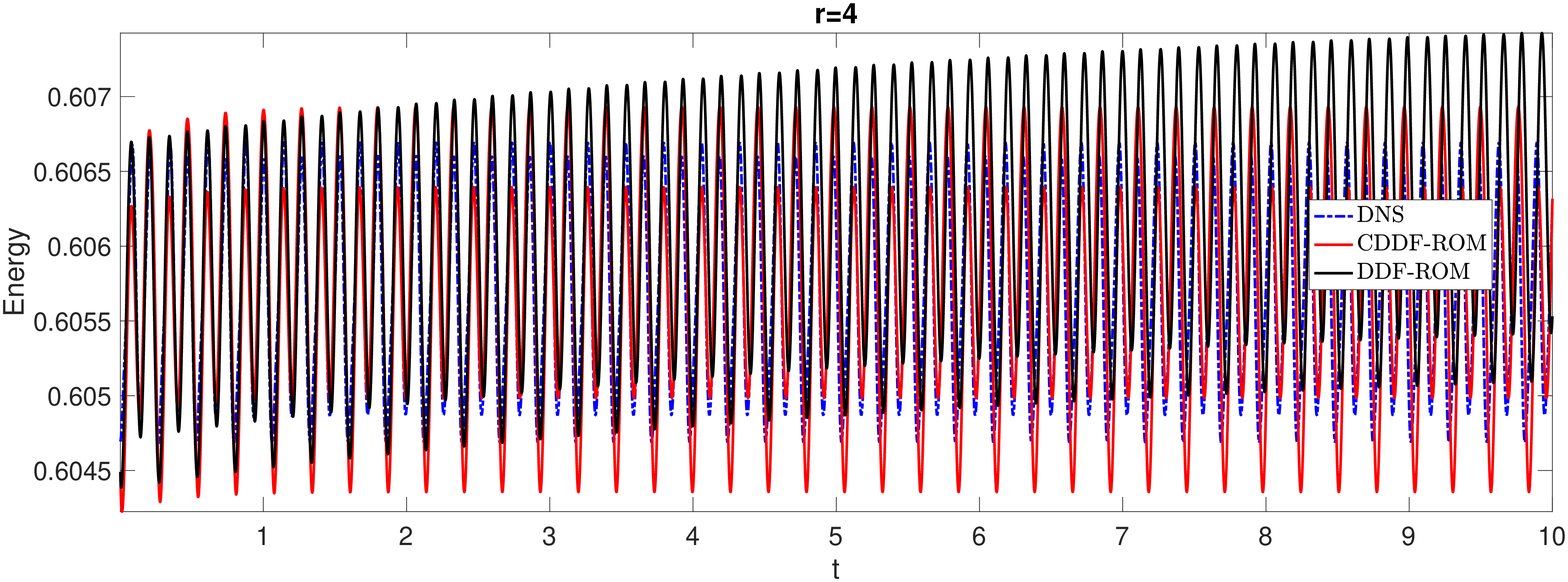}
     \caption{
     	Flow past a circular cylinder with $Re=1000$:
     	Plots of energy coefficients vs. time for DNS, DDF-ROM ($r=4$), and CDDF-ROM ($r=4$).
	\label{fig:Re1000r6}
	}
\end{figure}

Next, we investigate the case $r=4$.
This time, the minimum $m$ value that yields accurate results for CDDF-ROM and DDF-ROM is $m = r+3$.
For the CDDF-ROM, we use $tol= 0.42, \epsilon=2.5\times 10^{-3}$; for the DDF-ROM, we use $tol= 1\times 10^{-1}$.
In Fig.~\ref{fig:Re1000r6}, we plot the evolution of energy coefficients for CDDF-ROM, DDF-ROM, and FOM.
As in the $r=3$ case, the CDDF-ROM is significantly more accurate than the DDF-ROM.

\subsection{CDDF-ROM Cross-Validation: Predictive Investigation}
	\label{sec:cross-validation}

In this section, we perform a cross-validation of the CDDF-ROM. 
For comparison purposes, we also consider the DDF-ROM.
To this end, we test the two ROMs in settings that are different from the setting used to ``train'' the ROM closure model (i.e., $\btau$).
That is, we investigate the {\it predictive} capabilities of the CDDF-ROM. 

We consider two types of predictive investigations, both of which use fewer snapshots than those used in the previous sections: 

{\it (i) Equally Spaced Snapshots}: \ 
For this type of predictive investigation, in the constrained least squares problem \eqref{eqn:cddf-rom-least-squares} used to construct the CDDF-ROM's operators $\tilde{A}$ and $\tilde{B}$, instead of using ${\bf{a}}^{snap}(t_j)$ with $t_j, 1\le j\le M$ spanning an entire period, we used ${\bf{a}}^{snap}(t_{1+\ell(j-1)})$ with $t_{1+\ell(j-1)}, 1\le j\le \floor*{\frac{M+\ell-1}{\ell}}$ spanning just $\approx\frac{100}{\ell}\%$ of the entire period.
The parameter $M$ is the number of time steps required to complete a full period and takes the following values: 
For $Re=100, M=166$, for $Re=500, M=219$, and $Re=1000, M=134$.
 
{\it (ii) Unequally Spaced Snapshots}: \ 
This type of predictive investigation is similar to the one above.  
The only difference is that instead of choosing equally spaced snapshots, we use unequally spaced snapshots that are selected from the first part of the period.
Thus, the snapshots do not include information from the last part of the period, which makes this predictive investigation sonewhat more challenging than the previous one.

\subsubsection{Reynolds Number ${\bf Re=100}$}
In this section, we present results for $Re=100$.
We also consider $r=4$.

\begin{figure}[h]
 \includegraphics[width=.9\textwidth,viewport=0 25 1350 430]{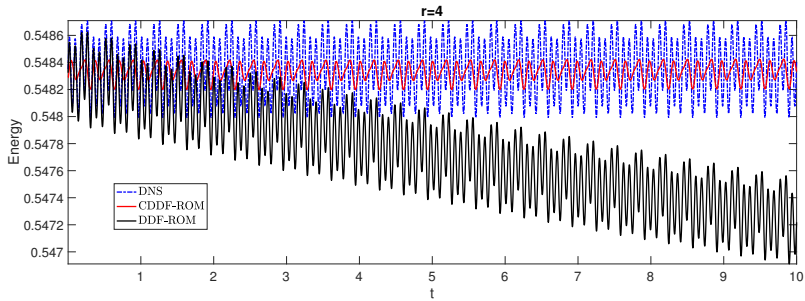}
     \caption{
     	CDDF-ROM cross-validation; 
	predictive investigation;
	equally spaced snapshots, $10.24\%$ data of an entire period;
     	flow past a circular cylinder with $Re=100$:
     	Plots of energy coefficients vs. time for DNS, DDF-ROM ($r=4$), and CDDF-ROM ($r=4$).
	\label{fig:prdc_N4d7m10}
	}
\end{figure}

First, we consider the case of equally spaced snapshots for $\ell=10$, which corresponds to $10.24\%$ data of an entire period.
We note that $m=r+3$ is the minimum $m$ value for which CDDF-ROM yields accurate results. 
For the CDDF-ROM, we use $tol=7\times 10^{-3}, \epsilon=1.94\times10^{-2}$; for the DDF-ROM, we use $tol= 7\times10^{-3}$.
In Fig.~\ref{fig:prdc_N4d7m10}, we plot the evolution of energy coefficients. 
Even for this drastic truncation, the CDDF-ROM performs very well; the DDF-ROM, on the other hand, is very inaccurate.

\begin{figure}[h]
 \includegraphics[width=.9\textwidth]{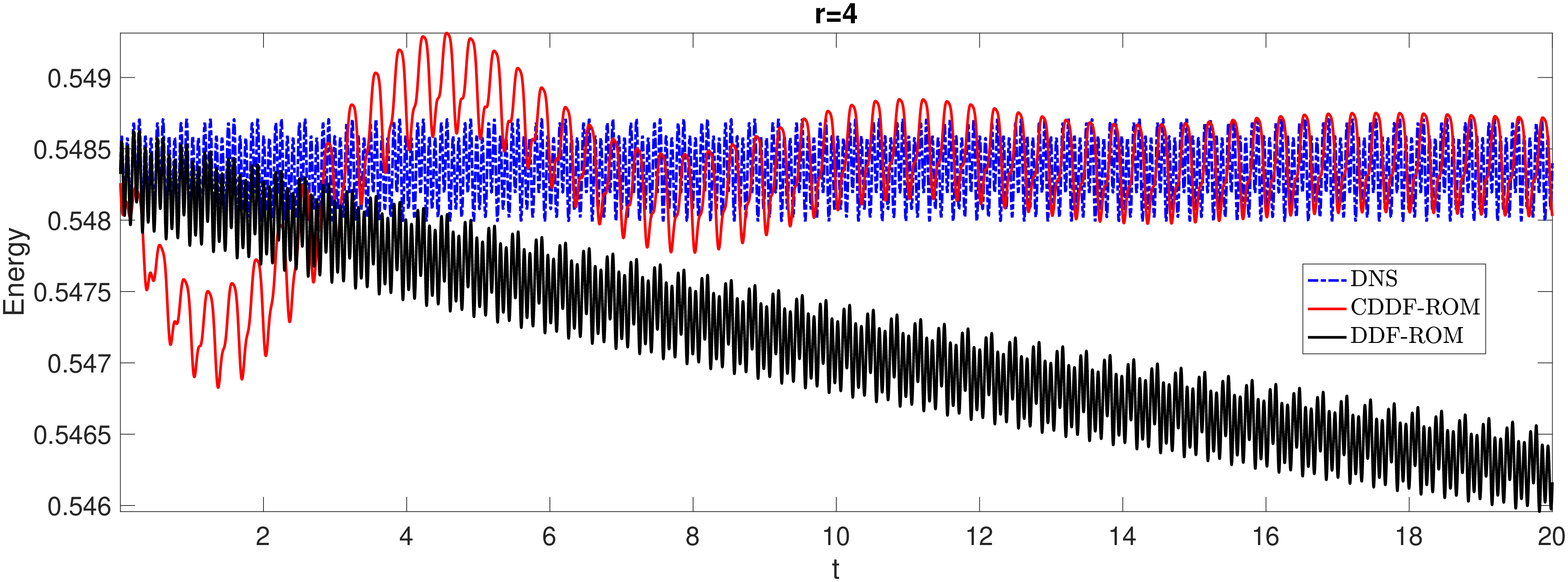}
     \caption{
     	CDDF-ROM cross-validation; 
	predictive investigation;
	unequally spaced snapshots, first $89\%$ of the entire period;
     	flow past a circular cylinder with $Re=100$:
     	Plots of energy coefficients vs. time for DNS, DDF-ROM ($r=4$), and CDDF-ROM ($r=4$).
	\label{fig3}
	}
\end{figure}

Next, we consider the case of unequally spaced snapshots collected from the first $89\%$ of the entire period.
This time, $m=r+1$ is the minimum $m$ value for which CDDF-ROM yields accurate results. 
For the CDDF-ROM, we use $tol=7\times 10^{-3}, \epsilon=10^{-10}$; for the DDF-ROM, we use $tol=7\times 10^{-3}$.
In Fig.~\ref{fig3}, we plot the evolution of energy coefficients.
The CDDF-ROM is again dramatically more accurate than DDF-ROM.

\subsubsection{Reynolds Number ${\bf Re=500}$}
In this section, we present results for $Re=500$.
Again, we consider $r=4$.

\begin{figure}[h]
	\includegraphics[width=.9\textwidth,viewport=0 25 1350 455, clip]{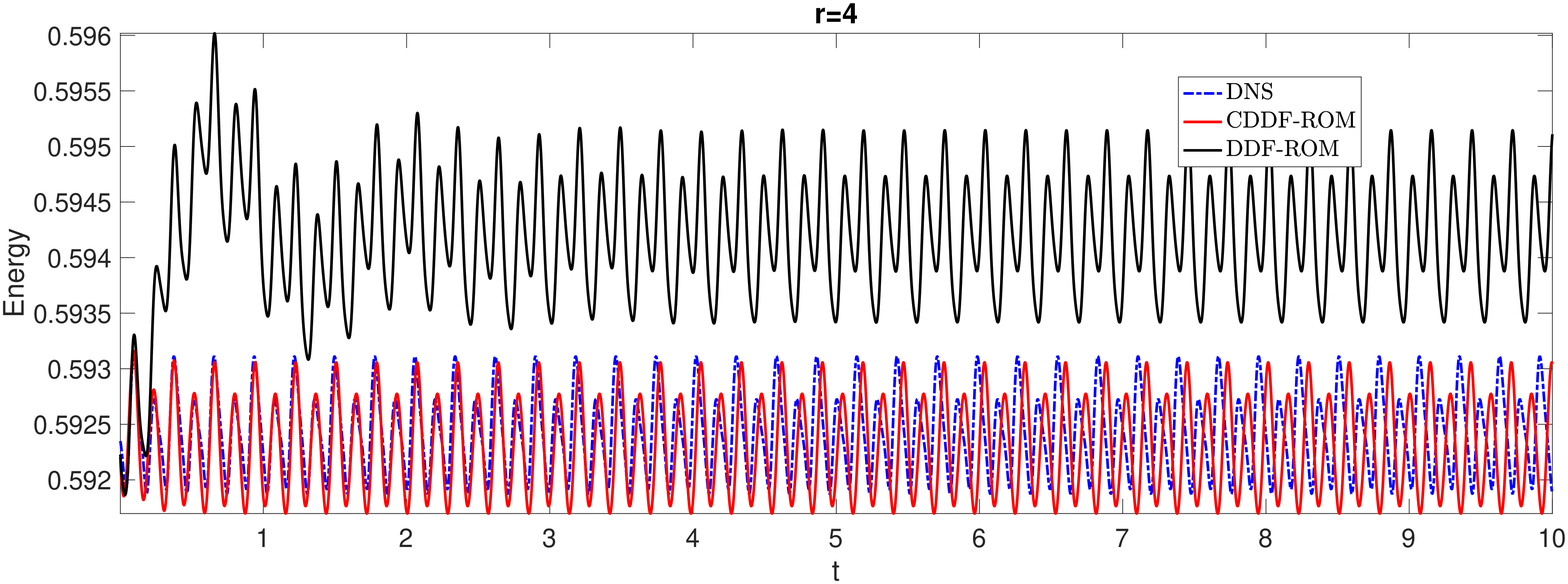}
     \caption{
     	CDDF-ROM cross-validation; 
	predictive investigation;
	equally spaced snapshots, $7\%$ data of an entire period;
     	flow past a circular cylinder with $Re=500$:
     	Plots of energy coefficients vs. time for DNS, DDF-ROM ($r=4$), and CDDF-ROM ($r=4$).
	\label{fig:predictive-equal-500}
	}
\end{figure}

First, we consider the case of equally spaced snapshots for $\ell=15$, which corresponds to $7\%$ data of an entire period.
We note that $m=r+1$ is the minimum $d$ value for which CDDF-ROM yields accurate results. 
For the CDDF-ROM, we use $tol= 10^{-3}, \epsilon=0$; for the DDF-ROM, we use $tol=3\times 10^{-3}$.
In Fig.~\ref{fig:predictive-equal-500}, we plot the evolution of energy coefficients. 
Even for this drastic truncation, the CDDF-ROM performs very well; the DDF-ROM, on the other hand, is very inaccurate.

\begin{figure}[h]
	\includegraphics[width=.9\textwidth,viewport=0 25 1350 455, clip]{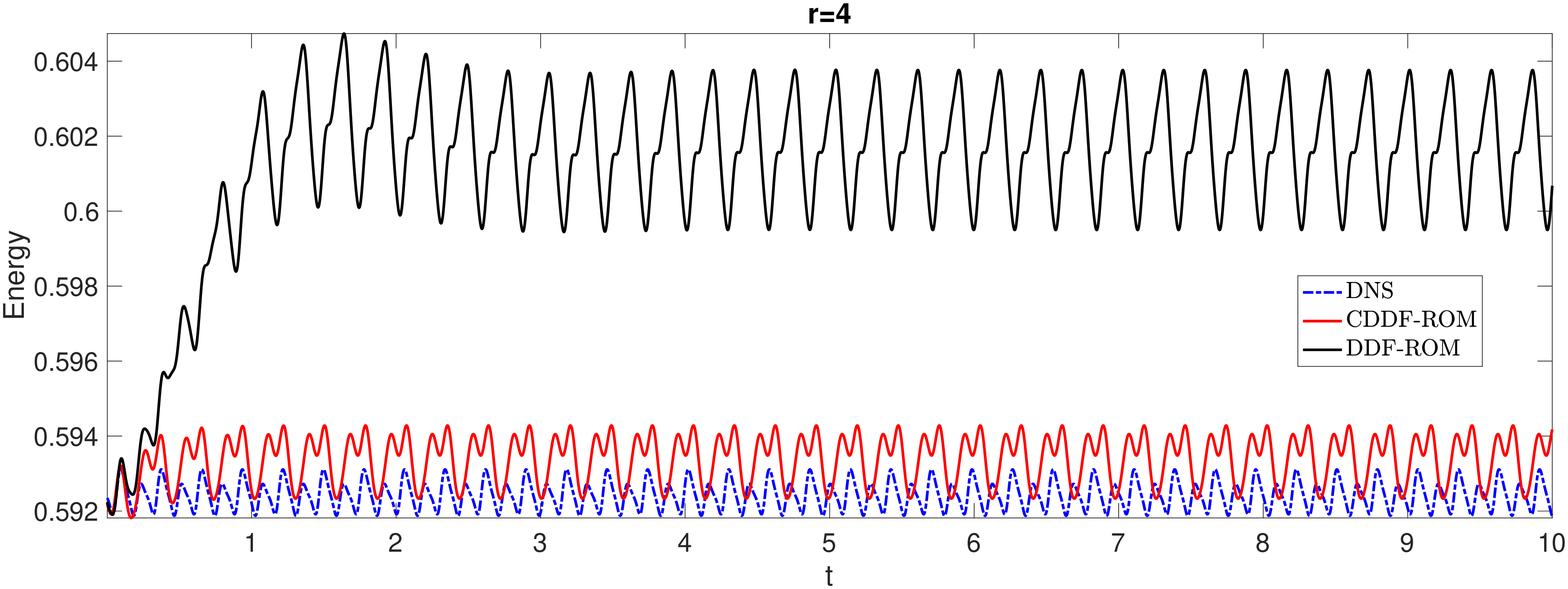}
     \caption{
     	CDDF-ROM cross-validation; 
	predictive investigation;
	unequally spaced snapshots, $50\%$ data of an entire period;
     	flow past a circular cylinder with $Re=500$:
     	Plots of energy coefficients vs. time for DNS, DDF-ROM ($r=4$), and CDDF-ROM ($r=4$).
	\label{fig:predictive-unequal-500}
	}
\end{figure}

Next, we consider the case of unequally spaced snapshots collected from the first $50\%$ of the entire period.
This time, $m=r+1$ is the minimum $d$ value for which CDDF-ROM yields accurate results. 
For the CDDF-ROM, we use $tol=3\times 10^{-4}, \epsilon=0$; for the DDF-ROM, we use $tol= 1\times 10^{-2}$.
In Fig.~\ref{fig:predictive-unequal-500}, we plot the evolution of energy coefficients.
The CDDF-ROM is again dramatically more accurate than DDF-ROM.

\subsubsection{Reynolds Number ${\bf Re=1000}$}
In this section, we present results for $Re=1000$.
This time, we consider $r=3$, which is the lowest $r$ value for which CDDF-ROM yields accurate results.

\begin{figure}[h]
	\includegraphics[width=.9\textwidth,viewport=0 25 1350 455, clip]{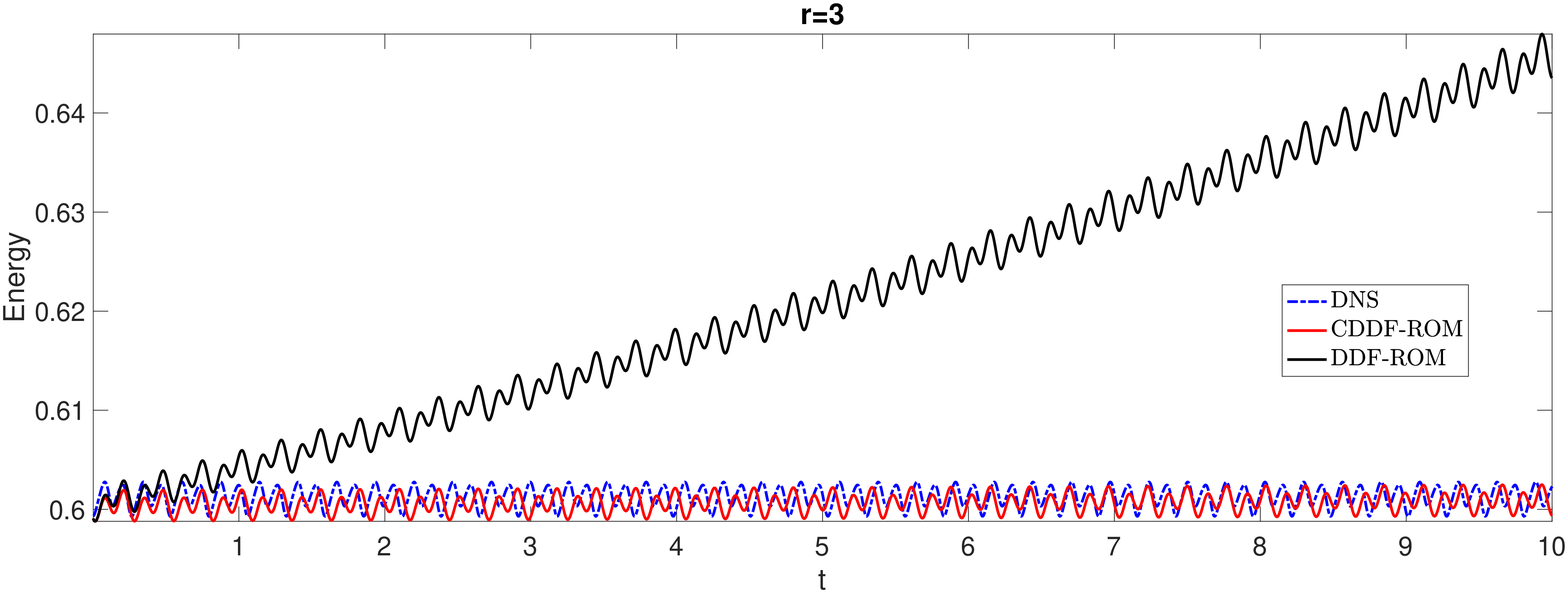}
     \caption{
     	CDDF-ROM cross-validation; 
	predictive investigation;
	equally spaced snapshots, $3\%$ data of an entire period;
     	flow past a circular cylinder with $Re=1000$:
     	Plots of energy coefficients vs. time for DNS, DDF-ROM ($r=3$), and CDDF-ROM ($r=3$).
	\label{fig:predictive-equal-1000}
	}
\end{figure}

First, we consider the case of equally spaced snapshots for $\ell=35$, which corresponds to $3\%$ data of an entire period.
We note that $m=r+1$ is the minimum $d$ value for which CDDF-ROM yields accurate results. 
For the CDDF-ROM, we use $tol=0.42, \epsilon=7\times 10^{-3}$; for the DDF-ROM, we use $tol= 5\times 10^{-3}$.
In Fig.~\ref{fig:predictive-equal-1000}, we plot the evolution of energy coefficients. 
Even for this drastic truncation, the CDDF-ROM performs very well; the DDF-ROM, on the other hand, is very inaccurate.

\begin{figure}[h]
	\includegraphics[width=.9\textwidth,viewport=0 25 1350 455, clip]{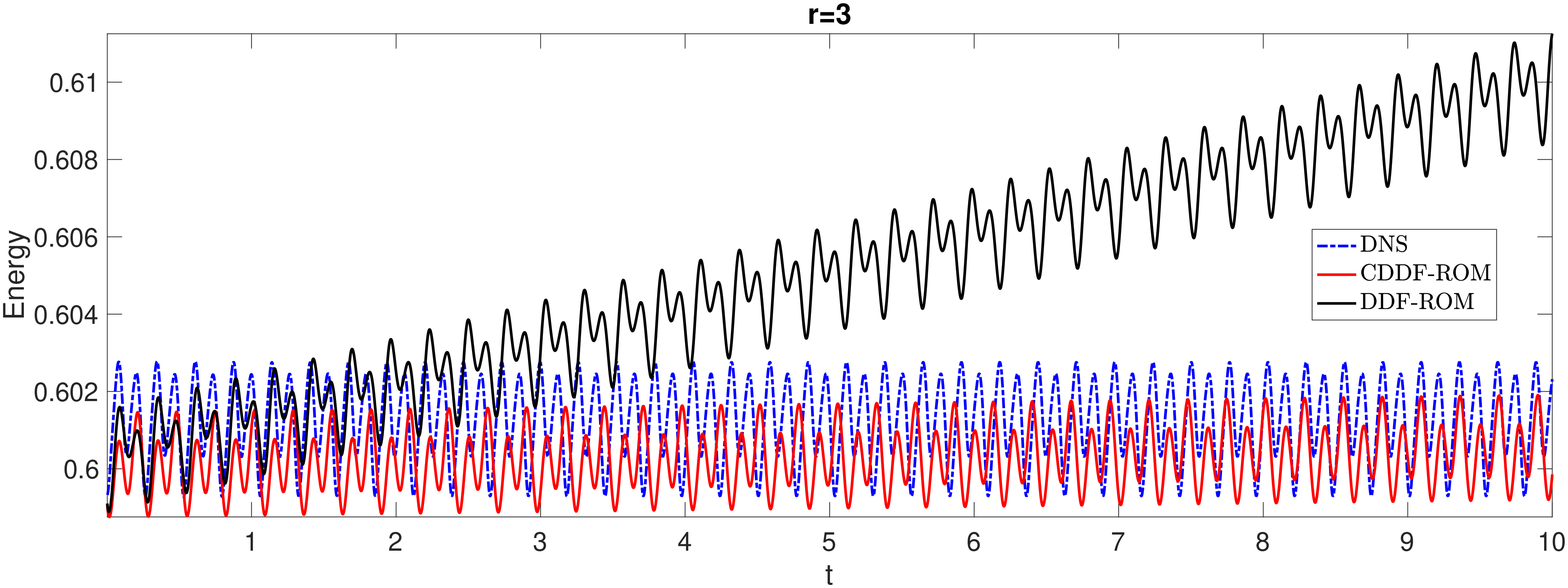}
     \caption{
     	CDDF-ROM cross-validation; 
	predictive investigation;
	unequally spaced snapshots, $50\%$ data of an entire period;
     	flow past a circular cylinder with $Re=1000$:
     	Plots of energy coefficients vs. time for DNS, DDF-ROM ($r=3$), and CDDF-ROM ($r=3$).
	\label{fig:predictive-unequal-1000}
	}
\end{figure}

Next, we consider the case of unequally spaced snapshots collected from the first $50\%$ of the entire period.
This time, $m=r+1$ is the minimum $d$ value for which CDDF-ROM yields accurate results. 
For the CDDF-ROM, we use $tol=0.42, \epsilon=2.5\times 10^{-2}$; for the DDF-ROM, we use $tol= 2.5$.
In Fig.~\ref{fig:predictive-unequal-1000}, we plot the evolution of energy coefficients.
The CDDF-ROM is again dramatically more accurate than DDF-ROM.

\subsubsection{Parameter Sensitivity}

To ensure a fair comparison between the DDF-ROM and CDDF-ROM, for each model we used optimal parameters $m$ (i.e., the parameter used in~\eqref{eqn:ddf-rom-efficiency} to ensure the CDDF-ROM's computational efficiency),  $tol$ (i.e., the tolerance value used in the truncated SVD method for the constrained least squares problem), and $\epsilon$ (i.e., the tolerance used to enforce the constraints $\tilde{A}_{ii} \leq -\epsilon$ in the constrained optimization problem solved with MATLAB toolbox lsqlin with interior-point method; see Section~\ref{sec:snapshot-rom-generation}). 
Those optimal parameters were selected to ensure that the ROM energy evolution matches best the FOM energy evolution.
(This, in turn, ensured the best ROM lift/drag evolutions.)

We note, however, that in our numerical investigation both the DDF-ROM and CDDF-ROM displayed a high sensitivity with respect to the parameters $m, tol$ and $\epsilon$. 
This is not surprising, since in Section 5.4 in~\cite{xie2018data} we showed that the DDF-ROM was sensitive with respect to parameters $m$ and $tol$.
The CDDF-ROM inherits this DDF-ROM's sensitivity and, in addition, is also sensitive with respect to $\epsilon$.
In Section~\ref{sec:conclusions}, we propose methods to alleviate the CDDF-ROM's sensitivity with respect to these parameters.

\section{Conclusions}
	\label{sec:conclusions}
	
In this paper, we have proposed a major improvement to the data-driven filtered ROM (DDF-ROM) introduced in~\cite{xie2018data} by adding physical constraints that increase the DDF-ROM's physical accuracy.
In the new physically-constrained DDF-ROM (CDDF-ROM), the constraints on the data-driven operators are mimicking the physical laws satisfied by the fluid flow equations (i.e., the nonlinear operator should conserve energy and the ROM closure term should be dissipative).
Thus, instead of using unconstrained data-driven modeling for the ROM closure problem (which is used in the DDF-ROM construction), we proposed using physically-constrained data-driven modeling.
We investigated the new CDDF-ROM in the numerical simulation of a 2D channel flow past a circular cylinder at Reynolds numbers $Re=100, Re=500$, and $Re=1000$.
To this end, we compared the CDDF-ROM with the DDF-ROM.
As a benchmark, we used the FOM data.
The CDDF-ROM was dramatically more accurate than the original DDF-ROM in predicting the evolution of the energy coefficients.
It was also more accurate in predicting the evolution of lift and drag coefficients, although the improvement was not as large.
Finally, we performed a {\it cross-validation} of the CDDF-ROM by testing it, together with the standard DDF-ROM, on data that was not used to train the ROM closure model; that is, we investigated the {\it predictive} capabilities of the CDDF-ROM.
To this end, we considered two types of cross-validations:
(i) with equally spaced snapshots, in which we considered as few as $3\%$ of the original set of snapshots; and 
(ii) with unequally spaced snapshots, in which we considered as few as $50\%$ of the first snapshots in the original set.
For both types of cross-validations, the CDDF-ROM dramatically outperformed the original DDF-ROM.
This numerical investigation in the reproductive and predictive regimes clearly showed that the physically-constrained CDDF-ROM is significantly more accurate than the DDF-ROM.  

We plan to investigate several research avenues: 
Probably the most important research direction is the further investigation of the CDDF-ROM's sensitivity with respect to the parameters used to solve the constrained least squares problem.
We plan to investigate alternative means of treating the ill-conditioning of the CDDF-ROM's constrained least squares problem (see, e.g., the approach used in~\cite{peherstorfer2016data}, where the authors combined trajectories of different initial conditions to decrease the ill-conditioning of their data-driven least squares problem).
Furthermore, we plan to investigate whether there is any connection between the ill-conditioning of the CDDF-ROM's constrained least squares problem and overfitting.
If such a connection exists, we plan to investigate methods that mitigate overfitting~\cite{brunton2016discovering,loiseau2018constrained,maulik2017neural}.
These alternative methods could eliminate altogether the need for the parameter $tol$ used in the CDDF-ROM's truncated SVD algorithm or yield numerical algorithms with lower parameter sensitivity.  

Another potential research avenue is the investigation of weaker constraints in the CDDF-ROM's constrained least squares problem.
For example, since the constraints~\eqref{eqn:atilde-constraint-componentwise-1}--\eqref{eqn:btilde-constraint-componentwise-3} are sufficient but not necessary conditions to satisfy the general constraint~\eqref{eqn:ddf-rom-constraint-atilde-btilde-separate}, we plan to investigate weaker constraints that satisfy the general constraint~\eqref{eqn:ddf-rom-constraint-atilde-btilde-separate}.
We will also investigate alternative models that impose the constraint~\eqref{eqn:ddf-rom-constraint-atilde-btilde} in a statistical sense, not pointwise, as is done in this paper.  
These alternative, weaker constraints might decrease the ill-conditioning of CDDF-ROM's constrained least squares problem, which in turn could yield numerical algorithms with lower parameter sensitivity.

\bibliographystyle{siamplain}
\bibliography{leo,traian}


\end{document}